# SPECTRAL GAPS IN WASSERSTEIN DISTANCES AND THE 2D STOCHASTIC NAVIER–STOKES EQUATIONS

By Martin Hairer[1] and Jonathan C. Mattingly[2]

*University of Warwick and Duke University*

We develop a general method to prove the existence of spectral gaps for Markov semigroups on Banach spaces. Unlike most previous work, the type of norm we consider for this analysis is neither a weighted supremum norm nor an $L^p$-type norm, but involves the derivative of the observable as well and hence can be seen as a type of 1-Wasserstein distance. This turns out to be a suitable approach for infinite-dimensional spaces where the usual Harris or Doeblin conditions, which are geared toward total variation convergence, often fail to hold. In the first part of this paper, we consider semigroups that have uniform behavior which one can view as the analog of Doeblin's condition. We then proceed to study situations where the behavior is not so uniform, but the system has a suitable Lyapunov structure, leading to a type of Harris condition. We finally show that the latter condition is satisfied by the two-dimensional stochastic Navier–Stokes equations, even in situations where the forcing is extremely degenerate. Using the convergence result, we show that the stochastic Navier–Stokes equations' invariant measures depend continuously on the viscosity and the structure of the forcing.

**1. Introduction.** This work is motivated by the study of the two-dimensional stochastic Navier–Stokes equations on the torus. However, the results and techniques are more general. The main abstract result of the paper gives a criterion guaranteeing that a Markov semigroup on a Banach space has a spectral gap in a particular 1-Wasserstein distance. (In the sequel, we will simply write Wasserstein for 1-Wasserstein.) To the best of our knowledge,

Received May 2006; revised May 2007.
[1]Supported in part by EPSRC Advanced Research Fellowship Grant EP/D071593/1.
[2]Supported in part by NSF PECASE Award DMS-04-49910 and an Alfred P. Sloan foundation fellowship.

these results are the first results providing a spectral gap in this, or any similar, setting. In turn, the existence of a spectral gap implies that the Markov semigroup possesses a unique, exponentially mixing invariant measure.

The results of this article rely on the existence of a "gradient estimate" introduced in [21] in the study of the degenerately forced Navier–Stokes equations on the two-dimensional torus. This estimate was used there in order to show that the corresponding Markov semigroup satisfies the "asymptotic strong Feller" property, also introduced in [21]. In this work, we show that gradient estimates of this form not only are useful to show uniqueness of the invariant measure, but are an essential ingredient to obtain the existence of a spectral gap for a large class of systems. In this introductory section, we concentrate on the two-dimensional stochastic Navier–Stokes equations on a torus to show how the main results can be applied. At the end of this section, we give an overview of the paper.

Recall that the Navier–Stokes equations describing the evolution of the velocity field $v(x,t)$ (with $x \in \mathbf{T}^2$) of a fluid under the influence of a body force $\bar{F}(x) + F(x,t)$ are given by

$$\text{(SNS)} \qquad \partial_t v = \nu \Delta v - (v \cdot \nabla)v - \nabla p + \bar{F} + F, \qquad \text{div} v = 0,$$

where the pressure $p(x,t)$ is determined by the algebraic condition $\text{div} v = 0$. We consider for $F$ a Gaussian stochastic forcing, that is, centered, white in time, colored in space and such that $\int \bar{F}(x) \, dx = \int F(x) \, dx = 0$. Since the gradient part of the forcing is canceled by the pressure term, we assume without loss of generality that $\text{div} \bar{F} = \text{div} F = 0$. More precisely, we assume that for $i, j \in \{1, 2\}$

$$\mathbf{E} F_i(x,t) F_j(x',t') = \delta(t-t') Q_{ij}(x-x'),$$

$$\sum_{i,j=1}^{2} \partial_{ij}^2 Q_{ij} = 0, \qquad \int Q_{ij}(x) \, dx = 0.$$

Although we are confident that our results are valid for $Q$ sufficiently smooth, we restrict ourselves to the case where $Q$ is a trigonometric polynomial, so that we can make use of the bounds obtained in [21, 39].

Instead of considering the velocity (SNS) directly, we will consider the equation for the vorticity $w = \nabla \wedge v = \partial_1 v_2 - \partial_2 v_1$. Note that $v$ is uniquely determined from $w$ (we will write $v = \mathcal{K}w$) through the conditions

$$w = \nabla \wedge v, \qquad \text{div} v = 0, \qquad \int v(x) \, dx = 0.$$

When written in terms of $w$, (SNS) is equivalent to

$$(1) \qquad \partial_t w = \nu \Delta w - (\mathcal{K}w) \cdot \nabla w + \bar{f} + f,$$

where we have defined $f = \nabla \wedge F$ and $\bar{f} = \nabla \wedge \bar{F}$. Note that since $f$ is translation invariant, one can write it as

$$f(x,t) = \text{Re} \sum_{k \in \mathbf{Z}^2 \setminus \{0\}} q_k e^{ikx} \xi_k(t),$$

where the $\xi_k$ are independent white noises and where $q_k = q_{-k}$. We can therefore identify the correlation function $Q$ with a vector $q$ in $\ell_+^2$, the set of square integrable sequences with positive entries. Denoting by $\mathcal{Z}$ the set of indices $k$ for which $q_k \neq 0$, we will make throughout this article the following assumptions:

ASSUMPTION 1. Only finitely many of the $q_k$'s are nonzero and $\bar{f}$ lies in the span of $\{e^{ikx} \mid q_k \neq 0\}$. Furthermore, $\mathcal{Z}$ generates $\mathbf{Z}^2$ and there exist $k, \ell \in \mathcal{Z}$ with $|k| \neq |\ell|$.

REMARK 1.1. The assumption that only a finite number of $q_k$ are nonzero is only a technical assumption reflecting a deficiency in [39]. All of the results of this article certainly hold if the first part of Assumption 1 is replaced by an appropriate decay property for the $q_k$. Note, for example, that in [21], Section 4.5, it is shown that there exists an $N_*$ such that if the range of $Q$ contains $\{e^{ikx} \mid |k| < N_*\}$, $\bar{f}$ is as in Assumption 1, and $\sum q_k^2 < \infty$, then all of the results of this paper hold. In particular, this allows infinitely many $q_k$ to be nonzero.

REMARK 1.2. Using the results in [3] one can remove the restriction that the forcing need consist of Fourier modes and replace it with the requirement that the forced functions span the Fourier modes required above. Since this is not the main point of this article, we do not elaborate further here.

REMARK 1.3. It is clear that the assumption that $\bar{f} \in \text{span}\{e^{ikx} \mid q_k \neq 0\}$ is far from optimal. The correct result likely places no restriction on $\bar{f}$ other than it be sufficiently smooth. This more delicate result requires an improved understanding of the control problem obtained by replacing the noise by controls. Some steps in this direction have been made [1, 2, 47], but the current results are not sufficient for our needs. Nonetheless, the present assumption on $\bar{f}$ seems reasonable from a modeling perspective where one would likely have some noise in all of the directions on which the body forces act.

We will consider (SNS) as an evolution equation in the subspace $\mathcal{H}$ of $H^1$ that consists of velocity fields $v$ with div$v = 0$ in the sense of distributions. Note that this is equivalent to $w \in \text{L}^2$. We make a slight abuse of notation and denote by $\mathcal{P}_t$ the transition probabilities for (SNS), as well as the corresponding Markov semigroup on $\mathcal{H}$, that is,

$$\mathcal{P}_t(v, A) = \mathbf{P}(v(t, \cdot) \in A \mid v(0, \cdot) = v),$$

for every Borel set $A \subset \mathcal{H}$, and

$$(\mathcal{P}_t \phi)(v) = \int_{\mathcal{H}} \phi(v') \mathcal{P}_t(v, dv'), \qquad (\mathcal{P}_t^* \mu)(A) = \int_{\mathcal{H}} \mathcal{P}_t(v, A) \mu(dv)$$

for every $\phi : \mathcal{H} \to \mathbf{R}$ and probability measure $\nu$ on $\mathcal{H}$. Analogously we define the projection $\mu \phi = \int \phi(x) \mu(dx)$. It was shown in [21] that Assumption 1 implies that (SNS) admits a unique invariant measure $\mu_\star$, that is, $\mu_\star$ is a probability measure on $\mathcal{H}$ such that $\mathcal{P}_t^* \mu_\star = \mu_\star$ for every $t \geq 0$.

This article is concerned with whether, for an arbitrary probability measure $\mu$, $\mathcal{P}_t^* \mu \to \mu_\star$ (as $t \to \infty$) and in what sense this convergence takes place. Note that (SNS) is not expected to have the strong Feller property, so that it is a fortiori not expected that $\mathcal{P}_t^* \mu \to \mu_\star$ in the total variation topology if $\mu$ and $\mu_\star$ are mutually singular. (See [9] for a general discussion of the strong Feller property in infinite dimensions and [21] for a discussion of its limitations in the present setting.)

In order to state the main result of the present article, we introduce the following norm on the space of smooth observables $\phi : \mathcal{H} \to \mathbf{R}$:

$$\|\phi\|_\eta = \sup_{v \in \mathcal{H}} e^{-\eta \|v\|^2} (|\phi(v)| + \|D\phi(v)\|).$$

Here, we denoted by $D\phi$ the Fréchet derivative of $\phi$. With this notation, we will show that the operator $\mathcal{P}_t$ has a spectral gap in the norm $\|\cdot\|_\eta$ in the following sense:

THEOREM 1.4. *Consider* (SNS) *in the range of parameters allowed by Assumption 1. For every $\eta$ small enough there exist constants $C$ and $\gamma$ such that*

$$\|\mathcal{P}_t \phi - \mu_\star \phi\|_\eta \leq C e^{-\gamma t} \|\phi\|_\eta,$$

*for every Fréchet differentiable function $\phi : \mathcal{H} \to \mathbf{R}$ and every $t \geq 0$.*

In [19] a similar operator-norm estimate on $\mathcal{P}_t$ was obtained in a weighted total variation norm ($\|\cdot\|_\eta$ without the $\|D\phi\|$ term) when the forcing was spatially rough and nondegenerate. Our setting is quite different. The spatially rough and nondegenerate forcing makes the analysis much closer to the finite-dimensional setting. It is not expected that such estimates hold in the total variation norm in the setting of this article. We should also remark that previous estimates, such as [20, 29, 37, 38], giving simply exponential convergence to equilibrium are weaker and the results in this article represent a significant and new extension of those results.

It is sometimes of interest to know that the structure functions of the solution to (SNS) converge to the structure functions determined by $\mu_\star$. This is not an immediate consequence of Theorem 1.4 because point evaluations of the velocity field are not continuous functions on $\mathcal{H}$. The smoothing properties of (SNS) nevertheless allow us to show the following result, which is an immediate consequence of Theorem 1.4 and Proposition 5.12 below.

THEOREM 1.5. *Consider* (SNS) *in the range of parameters allowed by Assumption 1. Let $n \geq 1$ and define the n-point structure functions by*

$$\mathcal{S}^n(x_1, \ldots, x_n) = \int v(x_1) \cdots v(x_n) \mu_\star(dv).$$

*Then, for every $\eta > 0$, there exist constants $C$ and $\gamma > 0$ such that, for every $v_0 \in \mathcal{H}$, one has the bound*

$$\sup_{x_1, \ldots, x_n} \left| \mathbf{E} \prod_{i=1}^n v(x_i, t) - \mathcal{S}^n(x_1, \ldots, x_n) \right| \leq C e^{\eta \|v_0\|^2 - \gamma t},$$

*for every $t > 1$. Here, $v(x, t)$ is the solution of* (SNS) *with initial condition $v_0$.*

The remainder of this article is structured as follows. In Section 2, we begin with an abstract discussion of our ideas in a setting with uniform estimates. In Section 3, we give the main theoretical results of the paper which combine the ideas from the first section with estimates stemming from an assumed Lyapunov structure. The convergence is measured in a distance in which paths are weighted by the Lyapunov function. We then turn in Section 5 to the specifics of the stochastic Navier–Stokes equation and show that it satisfies the general theorems from Section 3. In Section 5.3, we show that the Markov semigroup generated by (SNS) is strongly continuous on a suitable Banach space and that its generator has a spectral gap there. Then in Section 5.5, we demonstrate the power of the spectral gap estimates by giving a short proof that (SNS)'s invariant measures depend continuously on all the parameters of the equation.

**2. A simplified, uniform setting.** In this section, we illustrate many of the main ideas used throughout this article in a simplified setting. We consider the analogue of one of the simplest (and yet powerful) conditions for a Markov chain with transition probabilities $\mathcal{P}$ to have a unique invariant measure, namely Doeblin's condition[3]:

---

[3]Doeblin's original condition was the existence of a probability measure $\nu$ and a constant $\varepsilon > 0$ such that $\mathcal{P}(x, A) \geq \varepsilon$ whenever $\nu(A) > 1 - \varepsilon$; see [11, 12]. It turns out that, provided that the Markov chain is aperiodic and $\psi$-irreducible, this is equivalent to the (in general stronger) condition given here; see [40], Theorem 16.0.2.

THEOREM 2.1 (Doeblin). *Assume that there exist $\delta > 0$ and a probability measure $\nu$ such that $\mathcal{P}(x, \cdot) \geq \delta \nu$ for every $x$. Then, there exists a unique probability measure $\mu_\star$ such that $\mathcal{P}^* \mu_\star = \mu_\star$. Furthermore, one has $\|\mathcal{P}\phi - \mu_\star \phi\|_\infty \leq (1-\delta)\|\phi - \mu_\star \phi\|_\infty$ for every bounded measurable function $\phi$.*

A typical example of a semigroup for which Theorem 2.1 can be applied is given by a nondegenerate diffusion on a smooth compact manifold. Theorem 2.1 shows the fundamental mechanism for convergence to equilibrium in total variation norm. It is simple because the assumed estimates are extremely uniform. In this section we give a theorem guaranteeing convergence in a Wasserstein distance with assumptions analogous to Doeblin's condition.

A classical generalization of Doeblin's condition was made by Harris [22] who showed how to combine the existence of a Lyapunov function and a Doeblin-like estimate localized to a sufficiently large compact set to prove convergence to equilibrium. We will give a "Harris-like" version of our results in Section 3.

2.1. *Spectral gap under uniform estimates.* The aim of this section is to present a very simple condition that ensures that a Markov semigroup $\mathcal{P}_t$ on a Banach space $\mathcal{H}$ yields a contraction operator on the space of probability measures endowed with a Wasserstein distance. One can view it as a version of Doeblin's condition for the Wasserstein distance instead of the total variation distance. The main motivation for using a distance that metrizes the topology of weak convergence is that probability measures on infinite-dimensional spaces tend to be mutually singular, so that strong convergence is not expected to hold in general, even for ergodic systems.

The first assumption captures the regularizing effect of the Markov semigroup. While it does not imply that one function space is mapped into a more regular one as often occurs, it does say that at least gradients are decreased.

ASSUMPTION 2. *There exist constants $\alpha_1 \in (0,1)$, $C > 0$ and $T_1 > 0$ such that*

$$\|D\mathcal{P}_t \phi\|_\infty \leq C\|\phi\|_\infty + \alpha_1 \|D\phi\|_\infty, \tag{2}$$

*for every $t \geq T_1$ and every Fréchet differentiable function $\phi : \mathcal{H} \to \mathbf{R}$.*

REMARK 2.2. A typical way of checking (2) is to first show that for every $t \geq 0$, $\mathcal{P}_t$ is bounded as an operator on the space with norm $\|\phi\|_\infty + \|D\phi\|_\infty$. It then suffices to check that (2) holds with $\alpha_1 < 1$ for one particular time $t$ to deduce from the semigroup property that

$$\|D\mathcal{P}_t \phi\|_\infty \leq C(\|\phi\|_\infty + e^{-\gamma t}\|D\phi\|_\infty)$$

is valid with some $\gamma > 0$ for every $t \geq 0$.

REMARK 2.3. If the space $\mathcal{H}$ is actually a compact manifold, (2) together with the Arzelà–Ascoli theorem implies that the essential spectral radius of $\mathcal{P}_t$ (as an operator on the space with norm $\|\phi\|_\infty + \|D\phi\|_\infty$) is strictly smaller than 1. This is a well-known and often exploited fact[4] in the theory of dynamical systems. A bound like (2) is usually referred to as the Lasota–Yorke inequality [30, 31] or the Ionescu–Tulcea–Marinescu inequality [26] and is used to show the existence of absolutely continuous invariant measures when $\mathcal{P}_t$ is the transfer operator acting on densities. Usually, it is used with two different Hölder norms on the right-hand side. The present application with a Lipschitz norm and an infinity norm has a different flavor.

This bound alone is of course not enough in general to guarantee the uniqueness of the invariant measure. (Counterexamples with $\mathcal{H} = S^1$, the unit circle, can easily be constructed.) Furthermore, we are mainly interested in the case where $\mathcal{H}$ is not even locally compact.

In order to formulate our second assumption, we use the notation $\mathcal{C}(\mu_1 \mu_2)$ for the set of all measures $\Gamma$ on $\mathcal{H} \times \mathcal{H}$ such that $\Gamma(A \times \mathcal{H}) = \mu_1(A)$ and $\Gamma(\mathcal{H} \times A) = \mu_2(A)$ for every Borel set $A \subset \mathcal{H}$. Such a measure $\Gamma$ on the product space is referred to as a coupling of $\mu_1$ and $\mu_2$. We also denote by $\mathcal{P}_t^*$ the semigroup acting on probability measures which is dual to $\mathcal{P}_t$. With these notation, our second assumption, which is a form of uniform topological irreducibility, reads:

ASSUMPTION 3. For every $\delta > 0$, there exists a $T_2 = T_2(\delta)$ so that for any $t \geq T_2$ there exists an $a > 0$ so that

$$\sup_{\Gamma \in \mathcal{C}(\mathcal{P}_t^* \delta_x, \mathcal{P}_t^* \delta_y)} \Gamma\{(x', y') \in \mathcal{H} \times \mathcal{H} : \|x' - y'\| \leq \delta\} \geq a,$$

for every $x, y \in \mathcal{H}$.

REMARK 2.4. Recall that the total variation distance between probability measures can be characterized as one minus the supremum over all couplings of the mass of the diagonal. Therefore, if we set $\delta = 0$ in Assumption 3, we get $\|\mathcal{P}_t(x, \cdot) - \mathcal{P}_t(y, \cdot)\|_{\text{TV}} \leq 1 - a$ for every $x$ and $y$. By [40], Theorem 16.0.2, this is equivalent to the assumption of Theorem 2.1, so that the results in this section can be viewed as an analogue of Doeblin's theorem.

---

[4]It can be obtained as a corollary of the fact that the essential spectral radius of an operator $T$ can be characterized as the supremum over all $\lambda > 0$ such that there exists a singular sequence $\{f_n\}_{n \geq 0}$ with $\|f_n\| = 1$ and $\|T f_n\| \geq \lambda$ for every $n$. A slightly different proof can be found in [23] and is directly based on the study of the essential spectral radius by Nussbaum [42]. It is, however, very close in spirit to the much earlier paper [26].

To measure the convergence to equilibrium, we will use the following distance function on $\mathcal{H}$:

$$(3) \qquad d(x,y) = \min\{1, \delta^{-1}\|x-y\|\},$$

where $\delta$ is a small parameter to be adjusted later on. The distance (3) extends in a natural way to a Wasserstein distance between probability measures by

$$(4) \qquad d(\mu_1, \mu_2) = \sup_{\mathrm{Lip}_d(\phi)\leq 1}\left|\int \phi(x)\mu_1(dx) - \int \phi(x)\mu_2(dx)\right|,$$

where $\mathrm{Lip}_d(\phi)$ denotes the Lipschitz constant of $\phi$ in the metric $d$. By the Monge–Kantorovich duality [44, 49], the right-hand side of (4) is equivalent to

$$(5) \qquad d(\mu_1, \mu_2) = \inf_{\mu \in \mathcal{C}(\mu_1,\mu_2)} \iint d(x,y)\mu(dx,dy).$$

(Note that the infimum is actually achieved; see [50], Theorem 4.1.) With these notation, one has the following convergence result:

THEOREM 2.5. *Let $(\mathcal{P}_t)_{t\geq 0}$ be a Markov semigroup over a Banach space $\mathcal{H}$ satisfying Assumptions 2 and 3. Then, there exist constants $\delta > 0$, $\alpha < 1$ and $T > 0$ such that*

$$(6) \qquad d(\mathcal{P}_T^*\mu_1, \mathcal{P}_T^*\mu_2) \leq \alpha\, d(\mu_1,\mu_2),$$

*for every pair of probability measures $\mu_1$, $\mu_2$ on $\mathcal{H}$. In particular, $(\mathcal{P}_t)_{t\geq 0}$ has a unique invariant measure $\mu_\star$ and its transition probabilities converge exponentially fast to $\mu_\star$.*

PROOF. We will prove the general result by first proving it for delta measures, namely,

$$(7) \qquad d(\mathcal{P}_t^*\delta_x, \mathcal{P}_t^*\delta_y) \leq \alpha\, d(x,y)$$

for all $(x,y) \in \mathcal{H} \times \mathcal{H}$. Once this estimate is proven, we can finish the proof of the general case by the following argument.

The Monge–Kantorovich duality yields $Q \in \mathcal{C}(\mu_1,\mu_2)$ so that $d(\mu_1,\mu_2) = \int d(x,y)Q(dx,dy)$. To complete the proof observe that

$$\begin{aligned} d(\mathcal{P}_t^*\mu_1, \mathcal{P}_t^*\mu_2) &\leq \int d(\mathcal{P}_t^*\delta_x, \mathcal{P}_t^*\delta_y)Q(dx,dy) \\ &\leq \alpha \int d(x,y)Q(dx,dy) = \alpha\, d(\mu_1,\mu_2). \end{aligned}$$

Let us first show that (7) holds when $\|x - y\| \leq \delta$ for some appropriately chosen $\delta$. Note that by (4) this is equivalent to showing that

$$(8) \qquad |\mathcal{P}_t\phi(x) - \mathcal{P}_t\phi(y)| \leq \alpha\, d(x,y) \stackrel{\text{def}}{=} \alpha\delta^{-1}\|x-y\|$$

for all smooth $\phi$ with $\mathrm{Lip}_d(\phi) \leq 1$. Note that we can assume $\phi(0) = 0$, so that this implies that $\|D\phi\|_\infty \leq \delta^{-1}$ and $\|\phi\|_\infty \leq 1$. It follows from Assumption 2 that $\|D\mathcal{P}_t\phi\|_\infty \leq C + \alpha_1 \delta^{-1}$ for every $t \geq T_1$. Choosing $\delta = (1-\alpha_1)/(2C)$ and substituting in for $C$, we get $\|D\mathcal{P}_t\phi\|_\infty \leq \delta^{-1}(1+\alpha_1)/2$, so that (8) holds for $t \geq T_1$ and $\alpha \geq (1+\alpha_1)/2$.

Let us now turn to the case $\|x - y\| > \delta$. It follows from Assumption 3 that for every $t > T_2(\delta)$ there exists a positive constant $a$ so that for any $(x,y) \in \mathcal{H}^2$ there exists $\Gamma \in \mathcal{C}(\mathcal{P}_t^*\delta_x, \mathcal{P}_t^*\delta_y)$ such that $\Gamma(\Delta_\delta) > a > 0$, where

$$\Delta_\delta = \{(x',y') : \|x'-y'\| \leq \tfrac{1}{2}\delta\}.$$

Since $d(x',y') \leq \tfrac{1}{2}$ on $\Delta_\delta$ and $d(x',y') \leq 1$ on the complement, one has

$$\int d(x',y')\Gamma(dx',dy') \leq \frac{1}{2}\Gamma(\Delta_\delta) + 1 - \Gamma(\Delta_\delta) = 1 - \frac{1}{2}\Gamma(\Delta_\delta) \leq 1 - \frac{a}{2}.$$

Since we are in the setting $d(x,y) = 1$, this implies that when $\|x - y\| > \delta$,

$$|\mathcal{P}_t\phi(x) - \mathcal{P}_t\phi(y)| \leq \alpha\, d(x,y)$$

holds for $\alpha \geq 1 - \frac{a}{2}$ and $t \geq T_2(\delta)$.

Setting $\alpha = \max\{1-\frac{a}{2}, \frac{1}{2}(1+\alpha_1)\}$ and $T = \max\{T_1, T_2(\delta)\}$ completes the proof. $\square$

COROLLARY 2.6. *Let $\mathcal{P}_t$ be as in Theorem 2.5. Then, there exist constants $\alpha < 1$ and $T > 0$ such that*

$$(9) \quad \|\mathcal{P}_T\phi - \mu_\star\phi\|_{1,\infty} \leq \alpha\|\phi\|_{1,\infty}, \qquad \|\phi\|_{1,\infty} = \sup_{x \in \mathcal{H}}(|\phi(x)| + \|D\phi(x)\|),$$

*for every Fréchet differentiable function $\phi \colon \mathcal{H} \to \mathbf{R}$.*

PROOF. Define $d_1(x,y) = 1 \wedge \|x - y\|$. Since $d$ is equivalent to $d_1$, (6) still holds for arbitrary $\alpha$ (but with a different value for $T$) with $d$ replaced by $d_1$. The claim then follows from the Monge–Kantorovich duality, noting that $\mathrm{Lip}_{d_1}(\phi) \leq 2\|\phi\|_{1,\infty}$ and, for functions $\phi$ with $\int \phi(x)\mu_\star(dx) = 0$, $\|\phi\|_{1,\infty} \leq \mathrm{Lip}_{d_1}(\phi)$. $\square$

2.2. *A more pathwise perspective.* In [20, 37, 38], the authors advocated a pathwise point of view which explicitly constructed coupled versions of the process starting from two different initial conditions in such a way that the two coupled processes converged together exponentially fast. This point of view is very appealing as it conveys a lot of intuition; however, writing down the details can become a bit byzantine. Hence the authors prefer the arguments given in the preceding section for their succinctness and ease of verification. Nonetheless, the calculations of the present section provided the intuition which guided the previous section; and hence, we find it instructive to present them. As none of the rest of the paper uses any of the calculations from this section, we do not provide all of the details. Our goal is to show how the estimates from the previous section can be used to construct an explicit coupling in which the expectation of the distance between the trajectories starting from $x_0$ and $y_0$ converges to zero exponentially in time.

Fix a $t \geq \max(T_1, T_2)$, where $T_1$ and $T_2$ are the constants in Assumptions 2 and 3. Fix $\delta \geq 0$ as we did in the proof of Theorem 2.5. Now for $k = 0, 1, \ldots$ define the following sequence of stopping times:

$$r_k = \inf\{m \geq s_{k-1} : m \in \mathbf{N}, \|x_{mt} - y_{mt}\| \leq (1 - \alpha_1)\delta\},$$

$$s_k = \inf\{m \geq r_k : m \in \mathbf{N}, \|x_{mt} - y_{mt}\| > \delta\},$$

where $s_{-1} = 0$ by definition.

If $n \in [r_k, s_k)$, let the distribution of $(x_{t(n+1)}, y_{t(n+1)})$ be given by a coupling which minimizes the $\mathbf{E}\,d(x_{t(n+1)}, y_{t(n+1)})$. Hence choosing $\delta = (1 - \alpha_1)/(2C)$ and $\alpha \in ((1+\alpha_1)/2, 1)$ as in the paragraph below (8), Monge–Kantorovich duality gives that

(10) $$\mathbf{E}(d(x_{t(n+1)}, y_{t(n+1)})|(x_{tn}, y_{tn})) \leq \alpha \delta^{-1} \|x_{tn} - y_{tn}\|$$

provided $n \in [r_k, s_k)$. Given a random variable $X$ and events $A$ and $B$, for notational convenience, we define $\mathbf{E}(X; A) = \mathbf{E}(X \mathbf{1}_A)$ and $\mathbf{P}(A; B) = \mathbf{E}(\mathbf{1}_A \mathbf{1}_B)$, where $\mathbf{1}_A$ is the indicator function on the event $A$. Observe that if $\|x_0 - y_0\| \leq (1-\alpha_1)\delta$, then

(11) $$\mathbf{E}(\|x_{t(n+1)} - y_{t(n+1)}\|; s_1 > n+1|(x_{tn}, y_{tn})) \leq \alpha\|x_{tn} - y_{tn}\|\mathbf{1}_{s_1>n},$$

which implies that

$$\mathbf{E}(\|x_{t(n+1)} - y_{t(n+1)}\|; s_1 > n+1) \leq \alpha^{n+1}\|x_0 - y_0\|.$$

From this we see that as long as $x$ and $y$ stay in a $\delta$ ball of each other, they will converge toward each other exponentially in expectation.

Observe furthermore that

$$\delta\mathbf{P}(s_1 = n) = \delta\mathbf{P}(\|x_{tn} - y_{tn}\| > \delta; s_1 > n-1)$$

$$\leq \delta\mathbf{E}(d(x_{tn}, y_{tn}); s_1 > n-1)$$

$$= \delta\mathbf{E}(\mathbf{E}(d(x_{tn}, y_{tn}) \mid (x_{t(n-1)}, y_{t(n-1)})); s_1 > n-1)$$

$$\leq \alpha\mathbf{E}(\|x_{t(n-1)} - y_{t(n-1)}\|; s_1 > n-1) \leq \alpha^n\|x_0 - y_0\|,$$

where we used (10) to go from the third to the last line. Hence, assuming that $\|x_0 - y_0\| \leq (1-\alpha)\delta$,

$$(12) \qquad \mathbf{P}(s_1 > n) = 1 - \sum_{k=1}^{n-1} \mathbf{P}(s_1 = k) > 1 - \alpha,$$

so that $\mathbf{P}(s_1 < \infty) \leq \alpha < 1$. This shows that there is a positive chance that the two paths will indeed stay at distance less than $\delta$ from each other for all time.

All of the above calculations were predicated on the fact that $x_0$ and $y_0$ were initially less than $(1-\alpha)\delta$ apart. On the other hand, for $n \in [s_{k-1}, r_k)$, Assumption 3 guarantees that there exists a coupling for $(x_{t(n+1)}, y_{t(n+1)})$ so that

$$\mathbf{P}(\|x_{t(n+1)} - y_{t(n+1)}\| \leq (1-\alpha)\delta | (x_{tn}, y_{tn})) \geq a,$$

for some fixed $a > 0$. This shows that $\mathbf{P}(r_1 > n) \leq (1-a)^n$, so that the two paths will enter a $(1-\alpha)\delta$ ball of each other at a random time which has an exponentially decaying tail. We now sketch how to put these observations together more formally.

Let $d_1(x, y) = 1 \wedge \|x - y\|$ and define $\tau = \inf\{k : s_{k+1} = \infty\}$. We now combine the above estimates to sketch the proof of the exponential convergence to 0 of $\mathbf{E} d_1(x_{nt}, y_{nt})$. There are a few subtle issues arising from the fact that $\tau$ is not adapted to the natural filtration of the process, and we refer the interested reader to [20, 38, 43] for examples on how to circumvent these technicalities by a specific construction of $(x_{nt}, y_{nt})$. Since our goal is only to sketch the argument, we do not concern ourselves with these issues here.

Observe that for any $\beta \in (0, 1)$

$$\mathbf{E} d_1(x_{nt}, y_{nt}) \leq \mathbf{E}(d_1(x_{nt}, y_{nt}); r_\tau \leq n/2) + \mathbf{P}(\tau > \beta n)$$
$$+ \mathbf{P}(\tau < \beta n; r_\tau > n/2).$$

The first term decays exponentially fast in $n$ by (11), since the paths are guaranteed to be at distance less than $\delta$ on the time interval $[n/2, n]$. The second term is bounded by $\alpha^{\beta n}$ from the estimate $\mathbf{P}(s_1 < \infty) \leq \alpha$. Recall that the parameter $\beta$ is still free. Using the estimates from the preceding paragraph, it can be shown that for $\beta$ small enough the probability $\mathbf{P}(\tau < \beta n; r_\tau > n/2)$ has exponentially decaying tails since the random variable $r_{k+1} - r_k$ has exponentially decaying tails when restricted to the set where $s_k < \infty$.

**3. Spectral gap under a Lyapunov structure.** There are situations (the stochastic Navier–Stokes equations being a prime example) where it is not possible to verify Assumptions 2 and 3 in such a uniform way. The present section is an attempt to circumvent this by assuming that the system possesses a type of Lyapunov structure that compensates for the lack of uniformity of these estimates. The relationship between the results of the previous section and those of this section is analogous to the relationship between Doeblin's condition mentioned in the last section and Harris' conditions [16, 22, 40]. While the assumptions given in this section are heavily influenced by the known a priori bounds on the dynamics of the two-dimensional Navier–Stokes equations, we suspect the result will be useful more widely.

Throughout this section and the remainder of this article, we assume that we are given a random flow $\Phi_t$ on a Banach space $\mathcal{H}$. We will assume that the map $x \mapsto \Phi_t(\omega, x)$ is $\mathcal{C}^1$ for almost every element $\omega$ of the underlying probability space. We will denote by $D\Phi_t$ the Fréchet derivative of $\Phi_t(\omega, x)$ with respect to $x$.

Our first assumption is a strong type of Lyapunov structure on the flow:

ASSUMPTION 4. There exists a continuous function $V : \mathcal{H} \to [1, \infty)$ with the following properties:

1. There exist two strictly increasing continuous functions $V^*$ and $V_*$ from $[0, \infty) \to [1, \infty)$ so that

$$V_*(\|x\|) \leq V(x) \leq V^*(\|x\|) \tag{13}$$

   for all $x \in \mathcal{H}$ and such that $\lim_{a \to \infty} V_*(a) = \infty$.
2. There exist constants $C$ and $\kappa \geq 1$ such that

$$aV^*(a) \leq CV_*^\kappa(a), \tag{14}$$

   for every $a > 0$.
3. There exist a positive constant $C$, $r_0 < 1$, a decreasing function $\xi : [0, 1] \to [0, 1]$ with $\xi(1) < 1$ such that for every $h \in \mathcal{H}$ with $\|h\| = 1$

$$\mathbf{E}V^r(\Phi_t(x))(1 + \|D\Phi_t(x)h\|) \leq CV^{r\xi(t)}(x), \tag{15}$$

   for every $x \in \mathcal{H}$, every $r \in [r_0, 2\kappa]$ and every $t \in [0, 1]$.

REMARK 3.1. It follows from (15) and Jensen's inequality that there exists a constant $\tilde{C}$ such that

$$\mathbf{E}V^r(\Phi_t(x)) \leq \tilde{C}V^{r\xi(1)^{[t]}\xi(t-[t])}(x) \stackrel{\text{def}}{=} \tilde{C}V^{r\xi(t)}(x), \tag{16}$$

for every $t > 0$ and every $r \in [0, 2\kappa]$, where $[t]$ is the greatest integer smaller than $t$. In the last equality, we have extended the definition of $\xi$ to values of $t$ in $[0, \infty)$.

For $r \in (0,1]$, we introduce a family of distances $\rho_r$ on $\mathcal{H}$ by

$$\rho_r(x,y) = \inf_\gamma \int_0^1 V^r(\gamma(t)) \|\dot\gamma(t)\| \, dt,$$

where the infimum runs over all paths $\gamma$ such that $\gamma(0) = x$ and $\gamma(1) = y$.

In the interest of brevity, we will write $\rho$ for $\rho_1$. The main consequence of Assumption 4 used in this paper is that, via the distance function $\rho_r$, it also induces a kind of Lyapunov structure on the two-point dynamics:

LEMMA 3.2. *Assume that $\Phi_t$ is as above and that Assumption 4 holds. Then, for every $r \in [r_0, 1]$, there exist constants $\alpha \in (0,1)$ and $C, K > 0$ such that*

(17)
$$\mathbf{E}\rho_r(\Phi_t(x), \Phi_t(y)) \le C\rho_r(x,y),$$
$$\mathbf{E}\rho_r(\Phi_n(x), \Phi_n(y)) \le \alpha^n \rho_r(x,y) + K,$$

*for every $n \in \mathbf{N}$, every $t \in [0,1]$ and every pair $(x,y) \in \mathcal{H}^2$.*

PROOF. It suffices to show the second inequality in (17) for $n = 1$, since the other cases follow by iteration. Fix any $\epsilon > 0$ and fix a curve $\gamma$ connecting $x$ to $y$ such that

(18) $$\rho_r(x,y) \le \int_0^1 V^r(\gamma(t)) \|\dot\gamma(t)\| \, dt \le \rho_r(x,y) + \epsilon$$

and denote $\tilde\gamma(s) = \Phi_t(\gamma(s))$ for some $t \in [0,1]$. We then have

$$\mathbf{E}\rho_r(\Phi_t(x), \Phi_t(y)) \le \mathbf{E}\int_0^1 V^r(\tilde\gamma(s)) \|\dot{\tilde\gamma}(s)\| \, ds$$
$$\le \mathbf{E}\int_0^1 V^r(\tilde\gamma(s)) \|D\Phi_t(\gamma(s))\dot\gamma(s)\| \, ds$$
$$\le C\int_0^1 V^{r\xi(t)}(\gamma(s)) \|\dot\gamma(s)\| \, ds \le C\rho_r(x,y) + C\epsilon,$$

where the last inequality uses the fact that $\xi(t) \le 1$ by assumption. Since $\epsilon$ was arbitrary and $C$ independent of $\epsilon$ this yields the first bound in (17). Let now $R$ be sufficiently large so that $CV^{r\xi(1)}(x) \le \alpha V^r(x)$ for every $x$ with $|x| \ge R$. Such an $R$ exists since $V_\star$ tends to infinity. This yields

(19) $$\mathbf{E}\rho_r(\Phi_1(x), \Phi_1(y)) \le \alpha\rho_r(x,y) + CV^*(R) \int_0^1 \mathbf{1}_{B(R)}(\gamma(s)) \|\dot\gamma(s)\| \, ds + \epsilon,$$

where we denoted by $B(R)$ the ball of radius $R$ in $\mathcal{H}$ centered at the origin. Note now that

$$\int_0^1 \mathbf{1}_{B(R)}(\gamma(s)) \|\dot\gamma(s)\| \, ds \le 2R\left(\frac{V^*(R)}{V_*(0)}\right)^r + \epsilon,$$

since one could otherwise replace the corresponding piece of curve by a straight line and obtain a value which differed from $\rho_r(x, y)$ by more than $\epsilon$. Plugging this into (19) and again recalling the $\epsilon$ was arbitrary concludes the proof. $\square$

Our next assumption is a type of gradient inequality for the Markov semigroup $\mathcal{P}_t$ on $\mathcal{H}$ generated by the flow $\Phi_t$. In practice, this inequality can be verified if the system is hypoelliptic, in the sense of Hörmander (or effectively elliptic) and has suitable dissipative properties, but this is a hard task in general. (See [21] for a discussion of hypoellipticity and effective ellipticity in the setting of the 2D Navier–Stokes equations.)

ASSUMPTION 5. There exist a $C_1 > 0$ and $p \in [0, 1)$ so that for every $\alpha \in (0, 1)$ there exist positive $T(\alpha)$ and $C(\alpha)$ with

$$(20) \quad \|D\mathcal{P}_t \phi(x)\| \leq C_1 V^p(x)(C(\alpha)\sqrt{(\mathcal{P}_t|\phi|^2)(x)} + \alpha\sqrt{(\mathcal{P}_t\|D\phi\|^2)(x)}),$$

for every $x \in \mathcal{H}$ and $t \geq T(\alpha)$.

REMARK 3.3. While (20) is reminiscent of gradient estimates of the type considered in [4], there does not seem to be an obvious link between the two approaches. The main reason is that (20) is really a statement about the long-time behavior of $\mathcal{P}_t$ whereas the bounds in [4] are statements about the short-time behavior of $\mathcal{P}_t$.

Our final assumption is a relatively weak form of irreducibility:

ASSUMPTION 6. Given any $C > 0$, $r \in (0, 1)$ and $\delta > 0$, there exists a $T_0$ so that for any $T \geq T_0$ there exists an $a > 0$ so that

$$\inf_{|x|, |y| \leq C} \sup_{\Gamma \in \mathcal{C}(\mathcal{P}_T^* \delta_x, \mathcal{P}_T^* \delta_y)} \Gamma\{(x', y') \in \mathcal{H} \times \mathcal{H} : \rho_r(x', y') < \delta\} \geq a.$$

The main result of the present article is that under these conditions, one has uniform exponential convergence of the transition probabilities $\mathcal{P}_t(x, \cdot)$ to the (unique) invariant measure of the system:

THEOREM 3.4. Let $\Phi_t$ be a stochastic flow on a Banach space $\mathcal{H}$ which is almost surely $\mathcal{C}^1$ and satisfies Assumption 4. Denote by $\mathcal{P}_t$ the corresponding Markov semigroup and assume that it satisfies Assumptions 5 and 6. Then, there exist positive constants $C$ and $\gamma$ such that

$$(21) \qquad \rho(\mathcal{P}_t^* \mu_1, \mathcal{P}_t^* \mu_2) \leq Ce^{-\gamma t}\rho(\mu_1, \mu_2),$$

for every $t \geq 0$ and any two probability measures $\mu_1$ and $\mu_2$ on $\mathcal{H}$.

Since the space of probability measures $\mu$ on $\mathcal{H}$ such that $\rho(\mu, \delta_0) < \infty$ is complete for the topology induced by $\rho$ (see, e.g., [49]), (21) immediately yields:

COROLLARY 3.5. *Under the assumptions of Theorem 3.4, there exists a unique invariant probability measure $\mu_\star$ for $\mathcal{P}_t$.*

Before we turn to the proof of Theorem 3.4, we give a statement that is equivalent, but involves the semigroup acting on observables instead of the semigroup acting on measures. Since in this setting the semigroup $\mathcal{P}_t$ possesses an invariant measure $\mu_\star$, we can define the norm

$$(22) \qquad \|\phi\|_\rho = \sup_{x \neq y} \frac{|\phi(x) - \phi(y)|}{\rho(x,y)} + \left| \int_\mathcal{H} \phi(x) \mu_\star(dx) \right|.$$

An alternative definition of this norm is given in Lemma 4.2 in the next section.

Recall that we also make an abuse of notation by defining the projection operator $\mu_\star$ by $\mu_\star \phi = \int_\mathcal{H} \phi(y) \mu_\star(dy)$. With these notation, we have the following statement, which is the dual statement of Theorem 3.4:

THEOREM 3.6. *Let $\mathcal{P}_t$ be as in Theorem 3.4. Then, there exist constants $\gamma > 0$ and $C > 0$ such that*

$$\|\mathcal{P}_t \phi - \mu_\star \phi\|_\rho \leq C e^{-\gamma t} \|\phi - \mu_\star \phi\|_\rho,$$

*for every Fréchet differentiable function $\phi : \mathcal{H} \to \mathbf{R}$ and every $t > 0$.*

REMARK 3.7. This implies that the spectrum of $\mathcal{P}_t - \mu_\star$ as an operator on the space of Fréchet differentiable functions with finite $\|\cdot\|_\rho$-norm is contained in the disk of radius $e^{-\gamma t}$ around 0. Alternatively, $\mu_\star$ is an eigenvector for $\mathcal{P}_t^*$ with eigenvalue 1. All other eigenvectors have real part whose magnitude is at most $e^{-\gamma t}$. This is the gap referred to in the title of the article.

PROOF OF THEOREM 3.6. Since $\|\mathcal{P}_t \phi - \mu_\star \phi\|_\rho = \|\mathcal{P}_t(\phi - \mu_\star \phi)\|_\rho$, we can assume without loss of generality that $\mu_\star \phi = 0$. The claim then follows immediately from the fact that

$$|\mathcal{P}_T \phi(x) - \mathcal{P}_T \phi(y)| \leq \|\phi\|_\rho \rho(\mathcal{P}_T^* \delta_x, \mathcal{P}_T^* \delta_y) \leq C \|\phi\|_\rho e^{-\gamma t} \rho(x,y),$$

where the last inequality follows from Theorem 3.4. Dividing both sides by $\rho(x,y)$ and taking the supremum over $x$ and $y$ concludes the proof. □

3.1. *Proof of Theorem 3.4.* The proof of Theorem 3.4 is technically very simple but relies on a trick, which consists in considering instead of $\rho$ a distance $d$ which is equivalent to $\rho$ but behaves like a large constant times $\rho_r$ for nearby points and like a small constant times $\rho$ for points that are far apart.

More precisely, given three constants $\delta > 0$, $r \in [r_0, 1)$ and $\beta \in (0,1)$ to be determined later, we define

$$d(x,y) = \left(1 \wedge \frac{\rho_r(x,y)}{\delta}\right) + \beta \rho(x,y).$$

Note that $d$ is indeed equivalent to $\rho$ since $\rho_r \leq \rho$ and therefore

$$\beta \rho(x,y) \leq d(x,y) \leq (\delta^{-1} + \beta)\rho(x,y).$$

However, $d$ is much better than $\rho$ in capturing the geometry of the bounds available to us. This will allow us to proceed in a way similar to Section 2. This time, we will consider separately three cases. The first case, $\rho \geq K_\star$, $\rho_r \geq \delta$, will be treated by using the Lyapunov structure given by Lemma 3.2. The second case, $\rho_r < \delta$, will be treated by using the gradient estimate of Assumption 5. Finally, the last case, $\rho < K_\star$, $\rho_r \geq \delta$, will be treated using the irreducibility of Assumption 6. Lemma 3.9 is like the first part of the proof of Theorem 2.5 and Lemma 3.10 is like the second part. The first makes use of the local contraction guaranteed by Assumption 5. The second covers intermediate scales and uses Assumption 6 to ensure that the two points move close together some of the time to obtain a contraction estimate. Lemma 3.8 covers points far from the center of the space. Because of the weighting of the distance function by the Lyapunov function, there is contraction if the distant points simply move toward the center of the space.

The following three lemmas provide rigorous formulations of these claims.

LEMMA 3.8. *In the setting of Theorem 3.4, there exists a constant $K_\star$ such that for every $\delta > 0$, every $\beta \in (0,1)$ and every $r \in [r_0, 1)$ there exists a constant $\alpha_1 \in (0,1)$ such that*

$$\left.\begin{array}{r}\rho(x,y) \geq K_\star \\ \rho_r(x,y) \geq \delta\end{array}\right\} \implies \mathbf{E}\,d(\mathcal{P}_n^*\delta_x, \mathcal{P}_n^*\delta_y) \leq \alpha_1\,d(x,y)$$

*for all $n \in \mathbf{N}$.*

LEMMA 3.9. *In the setting of Theorem 3.4, for any $\alpha_2 \in (0,1)$ there exist $n_0 > 0$, $r \in [r_0, 1)$ and $\delta > 0$ so that*

$$\rho_r(x,y) < \delta \implies \mathbf{E}\,d(\mathcal{P}_n^*\delta_x, \mathcal{P}_n^*\delta_y) \leq \alpha_2\,d(x,y)$$

*for all $n \geq n_0$ and $\beta \in (0,1)$.*

LEMMA 3.10. *In the setting of Theorem 3.4, for any $K_\star$, $\delta > 0$, $r \in [r_0, 1)$ there exists an $n_1$ so that for any $n > n_1$ there is a $\beta \in (0,1)$ and an $\alpha_3 \in (0,1)$ so the following implication holds:*

$$\left.\begin{array}{l}\rho(x,y) < K_\star \\ \rho_r(x,y) \geq \delta\end{array}\right\} \implies \mathbf{E}\, d(\mathcal{P}_n^*\delta_x, \mathcal{P}_n^*\delta_y) \leq \alpha_3\, d(x,y).$$

It now suffices to show that the conditions of all three statements can be satisfied simultaneously in order to complete the proof of Theorem 3.4:

PROOF OF THEOREM 3.4. By the same argument as in the proof of Theorem 2.5, it suffices to prove that

$$(23) \qquad d(\mathcal{P}_t^*\delta_x, \mathcal{P}_t^*\delta_y) \leq \alpha\, d(x,y)$$

for all $(x,y) \in \mathcal{H} \times \mathcal{H}$.

We begin by fixing $K_\star$ as in Lemma 3.8. We then choose an arbitrary $\alpha_2 \in (0,1)$ and apply Lemma 3.9 which fixes $n_0 \geq 1$, $r \in [r_0, 1)$ and $\delta > 0$. With these in hand, we turn to Lemma 3.10 and fix an $N$ with $N \geq \max\{n_0, n_1\}$. This in turn fixes $\beta \in (0,1)$ and $\alpha_3 \in (0,1)$. Fixing $\beta$ sets the value of $\alpha_1$ in Lemma 3.8. Setting $\alpha = \max\{\alpha_1, \alpha_2, \alpha_3\} < 1$ completes the proof. □

We now turn to the proof of Lemmas 3.8–3.10.

PROOF OF LEMMA 3.8. It follows from Lemma 3.2 that there exist constants $\alpha_\star \in (0,1)$ and $K_\star > 0$ such that

$$\mathbf{E}\rho(\Phi_n(x), \Phi_n(y)) \leq \alpha_\star \rho(x,y),$$

for every $(x,y)$ such that $\rho(x,y) \geq K_\star$. Since $d(\mathcal{P}_n^*\delta_x, \mathcal{P}_n^*\delta_y) \leq \mathbf{E}\, d(\Phi_n(x), \Phi_n(y))$ we thus get the bound

$$d(\mathcal{P}_n^*\delta_x, \mathcal{P}_n^*\delta_y) \leq 1 + \alpha_\star \beta \rho(x,y).$$

On the other hand, $\rho_r(x,y) > \delta$ by assumption, so that

$$1 + \alpha_\star \beta \rho(x,y) = 1 - \alpha_\star + \alpha_\star\, d(x,y).$$

Since $d(x,y) \geq 1 + \beta K_\star$, this implies the claim with

$$\alpha_1 = \frac{1 + \alpha_\star \beta K_\star}{1 + \beta K_\star},$$

which is smaller than 1 (but close to it when $\beta$ is small) by construction. □

PROOF OF LEMMA 3.9. This lemma is the most delicate of the three in the sense that it does not follow from "soft" a priori estimates on the dynamic but requires to make use of the "hard," quantitative bound given by Assumption 5.

For the proof of this result, we use representation (4) for the distance $d$. Notice that we can assume without loss of generality that the test functions $\phi$ satisfy $\phi(0) = 0$ and are Fréchet differentiable, so that the condition $\text{Lip}_d(\phi) \leq 1$ together with (14) implies that

$$
\begin{aligned}
\|D\phi(x)\| &\leq (\delta^{-1} + \beta)V(x), \\
|\phi(x)| &\leq 1 + \beta\|x\|V^*(\|x\|) \leq 1 + \beta C V_*^\kappa(\|x\|).
\end{aligned}
\tag{24}
$$

Combining Assumption 5 with (24) and (16), we see that there exists a constant $C$ such that, for every $\alpha > 0$ there exists $C(\alpha)$ and $T_1(\alpha)$ such that

$$\|D\mathcal{P}_t\phi(x)\| \leq CV^{\kappa\xi(t)+p}(x)(C(\alpha) + \alpha\delta^{-1}),$$

for every $x \in \mathcal{H}$, every $t > T_1(\alpha)$ and every choice for $\delta$ and $\beta$ in $(0,1]$. Now fix an arbitrary value for $\alpha_3 \in (0,1)$ and pick $\alpha$ so that $\alpha C \leq \alpha_3/2$. By (15) there exists a $T(\alpha) \geq T_1(\alpha)$ so that $\kappa\xi(t) + p < 1$ for all $t \geq T(\alpha)$. At this point, we choose $r = \max\{r_0, \kappa\xi(T(\alpha)) + p\} < 1$. Using the above estimates produces

$$\|D\mathcal{P}_t\phi(x)\| \leq \delta^{-1}V^r(x)\left(\delta C(\alpha) + \frac{\alpha_3}{2}\right) \leq \alpha_3\delta^{-1}V^r(x),$$

where we choose $\delta$ sufficiently small in order to obtain the last inequality. Fixing any $\epsilon > 0$, let $\gamma \colon [0,1] \to \mathcal{H}$ be a curve connecting $x$ and $y$ as in (18) with $r = 1$. We have

$$
\begin{aligned}
|\mathcal{P}_t\phi(x) - \mathcal{P}_t\phi(y)| &= \left|\int_0^1 \langle D\mathcal{P}_t\phi(\gamma(s)), \dot{\gamma}(s)\rangle\, ds\right| \\
&\leq \alpha_3\delta^{-1}\int_0^1 V^r(\gamma(s))\|\dot{\gamma}(s)\|\, ds \\
&= \alpha_3\delta^{-1}\rho_r(x,y) + \epsilon\alpha_3\delta^{-1} \leq \alpha_3 d(x,y) + \epsilon\alpha_3\delta^{-1},
\end{aligned}
$$

where the last inequality uses the fact that we are in the case $\rho_r(x,y) \leq \delta$. Since $\epsilon$ was arbitrary, the proof is complete. □

In order to be able to prove Lemma 3.10 and thus conclude the proof of Theorem 3.4, it is essential to know that the region corresponding to the third case is a bounded subset of $\mathcal{H} \times \mathcal{H}$. This is given by the following result:

LEMMA 3.11. *Suppose that $V$ is as above and define, for some constants $\delta > 0$ and $K > 0$, the set*

$$\mathcal{C} = \{(x,y) \colon \rho_r(x,y) \geq \delta \text{ and } \rho(x,y) < K\}.$$

*Then, there exists an $R > 0$ such that $|x| \vee |y| \leq R$ for every $(x,y) \in \mathcal{C}$.*

PROOF. Note first that if $\delta/K > V_*^{r-1}(0)$, the set $\mathcal{C}$ is empty and there is nothing to prove.

We now show the contrapositive of the statement, that is, there exists $R > 0$ such that if $|x| > R$ and $\rho(x,y) < K$, then $\rho_r(x,y) < \delta$. Fixing any $\epsilon > 0$, let $\gamma$ denote a curve connecting $x$ to $y$ as in (18) with $r = 1$. Since $\rho(x,y) < K$ and $V \geq 1$, $\gamma$ never leaves the ball of radius $K + \epsilon$ around $x$. We thus have the bound

$$\rho_r(x,y) \leq \int_0^1 V^r(\gamma(s)) \|\dot\gamma(s)\| \, ds \leq \left(\sup_{z:\,|z-x|\leq K+\epsilon} V^{r-1}(z)\right)(\rho(x,y) + \epsilon)$$

$$= \left(\inf_{z:\,|x-z|\leq K+\epsilon} V(z)\right)^{r-1}(K + \epsilon).$$

Since $\epsilon$ was arbitrary and $V$ is continuous, the bound holds for $\epsilon = 0$. It follows from (13) that if one chooses $R = K + V_*^{-1}((\delta/K)^{1/(r-1)})$, one has

$$\left(\inf_{z:\,|x-z|\leq K} V(z)\right)^{r-1} \leq \frac{\delta}{K},$$

for every $x$ with $|x| > R$, which concludes the proof of the statement. □

With this fact secured, we are in the position to give the proof of Lemma 3.10.

PROOF OF LEMMA 3.10. Given $K_\star$, $\delta$ and $r \in (0,1)$, it follows from Lemma 3.11 that there exists a $C_*(K_\star, \delta, r)$ so that

$$\mathcal{C} \stackrel{\text{def}}{=} \{(x,y) : \rho_r(x,y) > \delta, \rho(x,y) \leq K_\star\} \subset \{(x,y) : \|x\|, \|y\| \leq C_*\}.$$

Hence by Assumption 6 for every $T$ large enough there exists a positive constant $a$ so that for any $(x_0, y_0) \in \mathcal{C}$ there exists a coupling $(x_T, y_T)$ of $(\Phi_T(x_0), \Phi_T(y_0))$ such that

$$\mathbf{P}(\rho_r(x_T, y_T) \leq \tfrac{1}{2}\delta) > a > 0.$$

Clearly $a$ is independent of the choice of $\beta$. Note now that there exists a constant $C$ such that, for every $z \in \mathcal{H}$,

$$\rho(z,0) \leq \int_0^1 V(sz)\|z\| \, ds \leq \|z\| V^*(\|z\|) \leq CV^\kappa(z).$$

Hence it follows from (15) that there exists a constant $C^*$ (also independent of $\beta$) such that $\mathbf{E}\rho(x_T, y_T) \leq \mathbf{E}\rho(x_T, 0) + \mathbf{E}\rho(y_T, 0) \leq C^*$ for all $(x_0, y_0) \in \mathcal{C}$.

As before, given a random variable $X$ and an event $A$, we define $\mathbf{E}[X; A] = \mathbf{E}[X\mathbf{1}_A]$. Now

$$\mathbf{E}\,d(x_T, y_T) = \mathbf{E}\left(1 \wedge \frac{\rho_r(x_T, y_T)}{\delta}; \rho_r(x_T, y_T) < \frac{1}{2}\delta\right)$$

$$+ \mathbf{E}\left(1 \wedge \frac{\rho_r(x_T, y_T)}{\delta}; \rho_r(x_T, y_T) \geq \frac{1}{2}\delta\right) + \beta \mathbf{E}\rho(x_T, y_T)$$
$$\leq \frac{1}{2}\left(1 - \mathbf{P}\left(\rho_r(x_T, y_T) \geq \frac{\delta}{2}\right)\right) + \mathbf{P}\left(\rho_r(x_T, y_T) \geq \frac{\delta}{2}\right)$$
$$+ \beta \mathbf{E}\rho(x_T, y_T)$$
$$\leq \frac{1}{2} + \frac{1}{2}\mathbf{P}(\rho_r(x_T, y_T) \geq \frac{1}{2}\delta) + \beta C^*$$
$$\leq \frac{1}{2} + \frac{1}{2}(1 - a) + \beta C^* = 1 - \frac{1}{2}a + \beta C^*.$$

By making $\beta$ small enough we can ensure that the right-hand side is less than 1. We denote this number by $\alpha_3$. Since $\rho_r(x,y) \geq \delta$ we know that $d(x,y) \geq 1$ and hence

$$\mathbf{E}\,d(x_T, y_T) \leq \alpha_3 \, d(x,y),$$

which is the quoted result. $\square$

**4. Quasi-equivalence of norms.** In the finite-dimensional setting where a Lyapunov function exists, it is natural to consider the norm on functions given by

$$\sup_x \frac{|\phi(x)|}{V(x)}. \tag{25}$$

(See, e.g., [40].) The norm on measures associated to it by duality is a weighted total variation norm. This norm can still be used in the infinite-dimensional setting *provided that the driving noise is sufficiently nondegenerate*; see, for example, [9] for a general theory and [19] for a recent ergodicity result on the stochastic two-dimensional Navier–Stokes equation.

In the present article, we are, however, interested in the situation where the driving noise is very degenerate. Indeed we assumed, for our main example of interest, that the driving noise is finite-dimensional, whereas the state space of our system is of course infinite-dimensional. In this setting, although it is possible to show that topological irreducibility holds, we do *not* expect the corresponding Markov process to be $\psi$-irreducible for any measure $\psi$. This is because, even though the system is formally hypoelliptic, we consider it very unlikely that it has the strong Feller property. It is indeed very simple to construct infinite-dimensional Ornstein–Uhlenbeck processes where the noise acts on every degree of freedom (so that the system is formally elliptic), but the system nevertheless lacks $\psi$-irreducibility. Therefore, the results from [40] are not applicable to the present situation and we do not expect to be able to get convergence results in the total variation norm.

It is therefore natural to look for a modification of (25) to the Wasserstein setting.

Motivated by these considerations, we introduce the following family of norms:

$$\|\phi\|_{V^r} = \sup_{x \in \mathcal{H}} \frac{|\phi(x)| + \|D\phi(x)\|}{V^r(x)}.$$

When we take $r = 1$, we will simply write $\|\phi\|_V$. The remainder of this section is devoted to showing that, modulo the semigroup $\mathcal{P}_t$, these norms can be considered to be equivalent to the norms $\|\cdot\|_{\rho_r}$ introduced in (22). Once this has been shown, we will have that Theorem 3.6 holds with the $\|\cdot\|_\rho$ norm replaced by the $\|\cdot\|_V$ norm defined above. This result is contained in Corollary 4.4. We begin by showing that the norm $\|\cdot\|_{\rho_r}$ is bounded from above and from below by the $\|\cdot\|_{V^{r'}}$ norm for a choice of $r'$ not necessarily equal to $r$.

PROPOSITION 4.1. *There exists a constant $C$ such that*

$$C^{-1}\|\phi\|_{V^{\kappa r}} \leq \|\phi\|_{\rho_r} \leq C\|\phi\|_{V^r},$$

*for every $r \in [0,1]$ and every Fréchet differentiable function $\phi$.*

Note first that:

LEMMA 4.2. *Recall the definition of $\|\cdot\|_{\rho_r}$ from (22) and let $\phi \colon \mathcal{H} \to \mathbf{R}$ be Fréchet differentiable. Then one also has*

$$(26) \qquad \|\phi\|_{\rho_r} = \sup_{x \in \mathcal{H}} \frac{\|D\phi(x)\|}{V^r(x)} + \int_{\mathcal{H}} \phi(x) \mu_\star(dx).$$

PROOF. Since

$$\lim_{\varepsilon \to 0} \sup_{y \colon \|y-x\| \leq \varepsilon} \frac{|\phi(x) - \phi(y)|}{\rho_r(x,y)} = \frac{\|D\phi(x)\|}{V^r(x)},$$

$\|\phi\|_{\rho_r}$ is greater than or equal to the right-hand side in (26). In order to prove the reverse inequality, we can assume without loss of generality that $\int \phi(x) \mu_\star(dx) = 0$ and $\|D\phi(x)\| \leq V^r(x)$ for all $x$. One then has

$$|\phi(x) - \phi(y)| = \int_0^1 \langle D\phi(\gamma(s)), \dot\gamma(s)\rangle \, ds \leq \int_0^1 V^r(\gamma(s))\|\dot\gamma(s)\| \, ds,$$

for any smooth path $\gamma$ connecting $x$ to $y$. Taking the infimum over all such $\gamma$ proves the claim. □

PROOF OF PROPOSITION 4.1. We start with the second inequality. It follows from Lemma 4.2 that it suffices to show that there exists $C > 0$ such that
$$\int \phi(x)\mu_\star(dx) \leq C\|\phi\|_V.$$
This follows immediately from the fact that $V$ is integrable against $\mu_\star$ by (15).

In order to show that the first inequality holds, fix $\phi$ with $\|\phi\|_{\rho_r} = 1$. One then has
$$|\phi(x) - \phi(0)| \leq \rho_r(x,0) \leq CV^{\kappa r}(x),$$
where the second inequality follows from (14). Furthermore, $\int \rho_r(x,0)\mu_\star(dx) \leq \int \rho(x,0)\mu_\star(dx) = C$. This yields
$$\left|\int \phi(x)\mu_\star(dx) - \phi(0)\right| \leq \int |\phi(x) - \phi(0)|\mu_\star(dx) \leq C,$$
so that $|\phi(0)| \leq C + \|\phi\|_{\rho_r} \leq C + 1$. Combining these bounds, we get
$$|\phi(x)| \leq |\phi(0)| + |\phi(x) - \phi(0)| \leq \tilde{C}V^{\kappa r}(x),$$
for some $\tilde{C} > 0$, which completes the proof. □

We now show that the semigroup $\mathcal{P}_t$ has the following contraction properties:

THEOREM 4.3. *There exist constants $C$ and $\gamma$ such that, for every $r \in [r_0, \kappa]$, every Fréchet differentiable function $\phi$ and every $t \geq 0$, one has the bounds*
$$\|\mathcal{P}_t\phi\|_{V^{r(t)}} \leq Ce^{\gamma t}\|\phi\|_{V^r}, \qquad \|\mathcal{P}_t\phi\|_{\rho_{r(t)}} \leq Ce^{\gamma t}\|\phi\|_{\rho_r},$$
*where $r(t) = \max\{\xi(t)r, r_0\}$.*

PROOF. It suffices to show the claims for $t \in [0,1]$ since the other cases follow by iteration. To begin with, we get bounds on the common term in both norms:
$$\|D\mathcal{P}_t\phi(x)\| \leq \mathbf{E}\|D\phi(\Phi_t(x))\|\|D\Phi_t(x)\|$$
$$\leq \left(\sup_{y \in \mathcal{H}} \frac{\|D\phi(y)\|}{V^r(y)}\right)CV^{r\xi(t)}(x),$$
where we made use of (15) in the last inequality. On the other hand, we have
$$\|\mathcal{P}_t\phi(x)\| \leq \left(\sup_{y \in \mathcal{H}} \frac{\|\phi(y)\|}{V^r(y)}\right)CV^{r\xi(t)}(x),$$

and, from the invariance of $\mu_\star$,

$$\int \mathcal{P}_t\phi(x)\mu_\star(dx) = \int \phi(x)\mu_\star(dx).$$

Combining these estimates proves the quoted results. □

COROLLARY 4.4. *There exist a time $T$ and a constant $C$ such that*

$$\|\mathcal{P}_T\phi\|_{V^r} \leq C\|\phi\|_{\rho_r},$$

*for every Fréchet differentiable function $\phi$ and every $r \in (\varepsilon/(1-\xi(1)), 1]$.*

PROOF. Let $r_n = \xi(1)^n \kappa r + \varepsilon(1-\xi(1)^n)/(1-\xi(1))$ as above. Then, we get from Theorem 4.3 and Proposition 4.1 that

$$\|\mathcal{P}_n\phi\|_{V^{r_n}} \leq C^n\|\phi\|_{V^{\kappa r}} \leq KC^n\|\phi\|_{\rho_r},$$

for some constants $C$ and $K$. Since we assume that $r > \varepsilon/(1-\xi(1)) = \lim_n r_n$, there exists $m$ such that $r_m \leq r$. The fact that $\|\phi\|_{V^r} \leq \|\phi\|_{V^{r_m}}$ completes the proof. □

An immediate consequence of Corollary 4.4 is the following result which states that Theorem 3.6 holds with $\|\cdot\|_\rho$ replaced by $\|\cdot\|_V$.

THEOREM 4.5. *Let $\mathcal{P}_t$ be as in Theorem 3.4. Then, there exist constants $\gamma > 0$ and $C > 0$ such that*

$$\|\mathcal{P}_t\phi - \mu_\star\phi\|_V \leq Ce^{-\gamma t}\|\phi - \mu_\star\phi\|_V,$$

*for every $\phi \in \mathcal{B}$ and every $t > 0$.*

**5. Application to the 2D stochastic Navier–Stokes equations.** We now apply the results of the previous sections to the two-dimensional Navier–Stokes equations on the torus $\mathbf{T}^2$, which is our main motivation for the present work. Recall that, in the vorticity formulation (1), these equations are given by

(27) $$dw = \nu\Delta w\,dt + B(\mathcal{K}w, w)\,dt + \bar{f}\,dt + Q\,dW(t),$$

$$w_0 \in \mathrm{L}^2(\mathbf{T}^2) \stackrel{\text{def}}{=} \mathcal{H},$$

where $B(u,v) = -(u\cdot\nabla)v$ is the usual Navier–Stokes nonlinearity, $W$ is a cylindrical Wiener process on $\mathcal{H}$, and $Q\colon\mathcal{H}\to\mathcal{H}$ is a positive self-adjoint finite rank operator commuting with translations. The viscosity $\nu > 0$ is arbitrary. We use the notation laid out in the Introduction. In particular, we denote by $e_k$, $k \in \mathbf{Z}^2$ the eigenfunctions of $\Delta$ and by $q_k$ the corresponding eigenvalues of $Q$. Unless indicated otherwise, we will assume that the constant component $\bar f$ of the body force and the coefficients $q_k$ satisfy Assumption 1.

It is well known (see, e.g., [8, 17]) that (27) has a unique solution under much weaker assumptions on the covariance operator $Q$. It is also well known that under similar conditions, (27) has an invariant measure $\mu_\star$. The uniqueness of this invariant measure is a much harder problem and has been a field of intense research over the past decade. Early results can be found in [9, 17, 36]. Until recently, the consensus that emerged in [5, 6, 14, 20, 28, 34, 37, 38] was that the uniqueness of the invariant measure can be obtained, provided that all the $q_k$ with $|k|^2 \leq N$ are nonzero, for some value $N \approx \sum q_k^2/\nu^3$. To the best of the author's knowledge, the only exception to this were the results of [15], that indicated that the invariant measure $\mu_\star$ should be unique provided that there exist $R > 0$ and $\alpha$ large enough such that all the $q_k$ with $|k| \geq R$ are bounded from above and from below by multiples of $|k|^{-\alpha}$. The uniqueness problem was eventually solved under Assumption 1 by the authors in the recent article [21]. This assumption is close to optimal since it only fails in situations where there exists a closed subspace $\widetilde{\mathcal{H}} \subset \mathcal{H}$ that is invariant for (27). It can then be shown that there always exists a unique ergodic invariant measure $\mu_\star$ for (27) such that $\mu_\star(\widetilde{\mathcal{H}}) = 1$.

We will show in this section that under Assumption 1, the random flow generated by the solutions of (27) satisfies the assumptions of Theorem 3.4 with $V(w) = \exp(\eta\|w\|^2)$ for a positive $\eta$ sufficiently small. We will then exhibit a Banach space of observables $\mathcal{B}$ which is such that the semigroup $\mathcal{P}_t$ generated by (27) extends to a *strongly continuous* semigroup of operators on $\mathcal{B}$. The results from Theorem 3.4 will then be shown to imply that the operator norm of $\mathcal{P}_t$ converges to 0, so that in particular its generator $\mathcal{L}$ has a spectral gap in the sense that there exists a constant $g > 0$ such that the spectrum of $\mathcal{L}$ is contained in $\{0\} \cup \{\mathrm{Re}\,\lambda \leq -g\}$. We conclude by showing first that $\mathcal{L}$ acts on cylindrical function as a second-order differential operator as one would expect and then that all the structure functions for (27) converge exponentially fast to their limit values.

5.1. *General Lyapunov structure.* We start with a result that we have found to be very useful when trying to check that (15) holds for a particular system.

LEMMA 5.1. *Let $U$ be a real-valued semimartingale*

$$dU(t,\omega) = F(t,\omega)\,dt + G(t,\omega)\,dB(t,\omega),$$

*where $B$ is a standard Brownian motion. Assume that there exist a process $Z$ and positive constants $b_1, b_2, b_3$, with $b_2 > b_3$, such that $F(t,\omega) \leq b_1 -$* $b_2 Z(t,\omega)$, $U(t,\omega) \leq Z(t,\omega)$ and $G(t,\omega)^2 \leq b_3 Z(t,\omega)$ almost surely. Then, the bound

$$\mathbf{E}\exp\left(U(t) + \frac{b_2 e^{-b_2 t/4}}{4}\int_0^t Z(s)\,ds\right) \leq \frac{b_2 \exp(2b_1/b_2)}{b_2 - b_3}\exp(U(0)e^{-(b_2/2)t})$$

holds for any $t \geq 0$.

PROOF. Fixing a time $t > 0$ and $a > 0$, set

$$Y(s) = \exp\left(\frac{b_2}{4}(s-t)\right)U(s) + \frac{b_2}{4}\int_0^s \exp\left(\frac{b_2}{4}(r-t)\right)Z(r)\,dr$$

and $M(s) = \int_0^s \exp(\frac{b_2}{4}(r-t))G(r,\omega)\,dB(r,\omega)$. Then

$$dY(s) = \exp\left(\frac{b_2}{4}(s-t)\right)\left(F(s,\omega) + \frac{b_2}{4}(U(s) + Z(s))\right)ds + dM(s).$$

If we restrict to $s \in [0,t]$, then we have that

$$Y(s) \leq Y(0) + \frac{4b_1}{b_2} - \frac{b_2}{2}\int_0^s \exp\left(\frac{b_2}{4}(r-t)\right)Z(r)\,dr + M(s).$$

Next observe that $Y(0) = \exp(-\frac{b_2}{4}t)U(0)$, $Y(t) \geq U(t) + \frac{b_2 e^{-b_2 t/4}}{4}\int_0^t Z(s)\,ds$ and

$$M(s) - \frac{b_2}{2}\int_0^s \exp\left(\frac{b_2}{4}(r-t)\right)Z(r)\,dr \leq M(s) - \frac{1}{2}\frac{b_2}{b_3}\langle M\rangle(s)$$

because $\exp(\frac{b_2}{4}(r-t))G(r^2) \leq \exp(\frac{b_2}{4}(r-t))Z(r)$ almost surely for $r \in [0,t]$. Since for continuous local martingales $M(t)$, one has the exponential martingale inequality

$$\mathbf{P}\left(\sup_s M(s) - \frac{\alpha}{2}\langle M\rangle(s) > \beta\right) \leq \exp(-\alpha\beta),$$

we have that

$$\mathbf{P}\left(U(t) + \frac{b_2 e^{-b_2 t/4}}{4}\int_0^t Z(s)\,ds - U(0)e^{-((b_2/2)t)} - \frac{2b_1}{b_2} > K\right) \leq \exp\left(-\frac{b_2}{b_3}K\right).$$

In order to conclude, it now suffices to use the fact that if $X$ is an arbitrary random variable and $a > 1$ is a constant such that $\mathbf{P}(X > K) \leq \exp(-aK)$ for every $K \geq 0$, then $\mathbf{E}\exp(X) \leq a/(a-1)$. □

5.2. *Verification of the assumptions of Theorem 3.4.* We first show that Lemma 5.1 indeed implies that:

PROPOSITION 5.2. *There exists $\eta_0$ such that, for every $\eta \in (0,\eta_0]$, the solutions to* (27) *satisfy Assumption 4 with $V(w) = \exp(\eta\|w\|^2)$.*

PROOF. It is clear that $V$ satisfies (13) and (14) so that it remains to show that (15) holds. Note that if we set $U(t) = \eta \|w\|^2$, we have from Itô's formula

$$dU(t) = \eta(\operatorname{tr} Q^2 + 2\langle w(t), \bar{f}\rangle - 2\nu\|w(t)\|_1^2)\,dt + 2\eta\|Qw(t)\|\,dB(t),$$

for some Brownian motion $B$. Here and in the sequel, we denote by $\|w\|$ the $L^2$-norm of $w$ and by $\|w\|_1 = \|\nabla w\|$ its $H^1$-norm. Since $\|w\|_1 \geq \|w\|$ and $2\langle w, \bar{f}\rangle \leq \nu^{-1}\|\bar{f}\|^2 + \nu\|w\|^2$, this shows that we are in the situation of Lemma 5.1 if we set $Z(t) = \eta\|w(t)\|_1^2$ and

$$b_1 = \eta\operatorname{tr} Q^2 + \frac{\|\bar{f}\|^2}{\nu}, \qquad b_2 = \nu, \qquad b_3 = 4\eta\|Q\|.$$

In particular, this shows that, for every $\eta < \nu/(4\|Q\|)$, there exists a constant $C$ such that, for every $t \in [0,1]$,

$$(28) \quad \mathbf{E}\exp\left(\eta\|w(t)\|^2 + \frac{\nu\eta e^{-\nu/2}}{2}\int_0^t \|w(s)\|_1^2\,ds\right) \leq C\exp(\eta\|w(0)\|^2 e^{-(\nu t/2)}).$$

On the other hand, we know from Lemma A.1 that, for every $\kappa > 0$, there exists a constant $C$ such that

$$\|D\Phi_t(w_0)\| \leq C\exp\left(\kappa\int_0^t \|w(s)\|_1^2\,ds\right) \qquad \forall t \in [0,1],$$

holds almost surely for every $w \in \mathcal{H}$. Combining this with (28) shows that (15) holds with $\xi(t) = e^{-\nu t/2}$ for arbitrarily small values of $r_0$. □

Recall now that the following "gradient estimate" is the main technical result of [21]:

PROPOSITION 5.3. *For every $\eta > 0$ and every $\alpha > 0$, there exist constants $C_{\eta,\alpha}$ such that, for every Fréchet differentiable function $\phi$ from $\mathcal{H}$ to $\mathbf{R}$, one has the bound*

$$\|D\mathcal{P}_n\phi(w)\| \leq \exp(\eta\|w\|^2)(C_{\eta,\alpha}\sqrt{(\mathcal{P}_n|\phi|^2)(w)} + \alpha^n\sqrt{(\mathcal{P}_n\|D\phi\|^2)(w)}),$$

*for every $w \in \mathcal{H}$ and $n \in \mathbf{N}$.*

REMARK 5.4. The works [21, 39] made the assumption $\bar{f} = 0$. However, the arguments presented there work without any modification under the assumption that $\bar{f} \in \operatorname{range} Q$. Note, for example, that Girsanov's formula implies that the transition probabilities for (SNS) with $\bar{f} = 0$ are equivalent to the transition probabilities with $\bar{f} \in \operatorname{range} Q$. In particular, this means that the proof of weak irreducibility from [21] carries over to the setting of this paper.

Proposition 5.3 immediately implies that Assumption 5 is satisfied for every choice of $\eta$, so that it remains to verify Assumption 6. This, however, follows immediately from [13], Lemma 3.1, and Remark 5.4 above. As a consequence, we have just shown that:

THEOREM 5.5. *If Assumption 1 holds, there exists $\eta_0 > 0$ such that, for every $\eta \leq \eta_0$, the stochastic flow solving* (27) *satisfies the assumptions of Theorem 3.4 with $V(w) = \exp(\eta \|w\|^2)$. Hence, the conclusions of Theorems 3.4, 3.6 and 4.5 hold.*

5.3. *Spectral gap for the generator.* In this section, we show that it is possible to extend the Markov semigroup $\mathcal{P}_t$ generated by solutions to (27) to some Banach space of observables $\mathcal{B}$ in such a way that:

1. The semigroup $\mathcal{P}_t$ is strongly continuous on $\mathcal{B}$.
2. There exists $g > 0$ such that $\sigma(\mathcal{P}_t) \setminus \{1\}$ is included in the disk of radius $e^{-gt}$ for every $t > 0$. Here, $\sigma(\mathcal{P}_t)$ denotes the spectrum of $\mathcal{P}_t$ viewed as a bounded operator on $\mathcal{B}$.

REMARK 5.6. It follows from standard semigroup theory that the above statements imply that $\mathcal{P}_t$ possesses a generator $\mathcal{L}$ densely defined on $\mathcal{B}$ (see, e.g., [10], Theorem 1.7) and that there exists $g > 0$ such that $\text{Re}(\lambda) \leq -g$ for every $\lambda \in \sigma(\mathcal{L}) \setminus \{0\}$ (see, e.g., [10], Theorem 2.16).

Before we give the precise statement of our results, let us turn to the construction of the Banach space $\mathcal{B}$. Given a Hilbert space $\mathcal{H}$, we define $\mathcal{C}_0^\infty(\mathcal{H})$ by

$$\mathcal{C}_0^\infty(\mathcal{H}) = \{\phi \circ \Pi | \Pi : \mathcal{H} \to \mathbf{R}^n \text{ linear}, \phi \in \mathcal{C}_0^\infty(\mathbf{R}^n)\}.$$

Note in particular that elements of $\mathcal{C}_0^\infty(\mathcal{H})$ are Fréchet differentiable of all orders. Given $\eta > 0$, define $\mathcal{B}_\eta$ as the closure of $\mathcal{C}_0^\infty(\mathcal{H})$ under the norm

$$(29) \qquad \|\phi\|_\eta = \sup_{w \in \mathcal{H}} \exp(-\eta \|w\|^2)(|\phi(w)| + \|D\phi(w)\|).$$

We also denote by $\tilde{\mathcal{B}}_\eta$ the closure under this norm of the space of all Fréchet differentiable functions $\phi$ such that $\|\phi\|_\eta$ is finite.

REMARK 5.7. The space $\mathcal{B}_\eta$ is much smaller than $\tilde{\mathcal{B}}_\eta$. In particular, elements of $\mathcal{B}_\eta$ are continuous when $\mathcal{H}$ is equipped with the topology of weak convergence, so that $w \mapsto \|w\|^2$ does *not* belong to $\mathcal{B}_\eta$, even though it obviously belongs to $\tilde{\mathcal{B}}_\eta$. However, $w \mapsto \|Kw\|^2$ does belong to $\mathcal{B}_\eta$, provided that $K : \mathcal{H} \to \mathcal{H}$ is a compact operator.

REMARK 5.8. The fact that the vorticity belongs to $\mathcal{H} = \mathrm{L}^2$ does not ensure that the corresponding velocity field is continuous. Therefore, point evaluations of the velocity field do not belong to $\mathcal{B}_\eta$. This fact can, however, be dealt with and we will do so in Section 5.4.

REMARK 5.9. Given an orthonormal basis $\{e_n\}$ of $\mathcal{H}$, one could have restricted oneself to the set of all functions of the type $w \mapsto \phi(\langle w, e_1 \rangle, \ldots, \langle w, e_n \rangle)$ with $\phi \in \mathcal{C}_0^\infty(\mathbf{R}^n)$. It is easy to check that the closure of this set under the norm (29) is again equal to $\mathcal{B}_\eta$, independently of the choice of basis.

As a consequence of this, it is a straightforward exercise to check that polynomials in $\langle w, e_n \rangle$ with rational coefficients form a dense subset of $\mathcal{B}_\eta$, so that it is a separable Banach space.

The first result of this section is the following:

THEOREM 5.10. *For $\eta$ sufficiently small, $\mathcal{P}_t$ extends to a $C_0$-semigroup on $\mathcal{B}_\eta$.*

PROOF. Define $\Pi_n$ as the orthogonal projection in $\mathcal{H}$ onto the first $n$ Fourier modes. The proof of this result is broken into two distinct steps as follows:

1. The semigroup $\mathcal{P}_t$ extends to a semigroup of bounded operators on $\mathcal{B}_\eta$ that is uniformly bounded as $t \to 0$.
2. One has $\|\mathcal{P}_t \phi - \phi\|_\eta \to 0$ as $t \to 0$ for a dense subset of elements of $\mathcal{B}_\eta$.

Note first that it follows from the a priori bounds of Lemma A.1 that if $\phi: \mathcal{H} \to \mathbf{R}$ is a Fréchet differentiable function such that $\|\phi\|_\eta < \infty$, then $\mathcal{P}_t \phi$ is again Fréchet differentiable and there exist constants $C_t$ that remain bounded as $t \to 0$ such that

$$\|\mathcal{P}_t \phi\|_\eta \leq C_t \|\phi\|_\eta,$$

provided that $\eta$ is sufficiently small. This shows that $\mathcal{P}_t$ can be extended to a semigroup on $\tilde{\mathcal{B}}_\eta$ which is uniformly bounded as $t \to 0$.

Since the norm on $\tilde{\mathcal{B}}_\eta$ is the same as on $\mathcal{B}_\eta$, the first claim follows if we can show that $\mathcal{P}_t$ maps $\tilde{\mathcal{B}}_\eta$ into itself. For an arbitrary function $\phi \in \mathcal{C}_0^\infty(\mathcal{H})$, we will show that

$$\lim_{n \to \infty} \|\mathcal{P}_t \phi - (\mathcal{P}_t \phi) \circ \Pi_n\|_\eta = 0, \tag{30}$$

where $\Pi_n$ denotes the orthogonal projection in $\mathcal{H}$ onto the Fourier modes with $|k| \leq n$. This is sufficient since it follows from the a priori bounds (46), (43), (40) and (42) that the function $(\mathcal{P}_t \phi) \circ \Pi_n$ is twice Fréchet differentiable and that, together with its derivative, it grows more slowly than $\exp(\eta \|x\|^2)$ at infinity, so that it belongs to $\mathcal{B}_\eta$.

Fix a generic element $w \in \mathcal{H}$ and a natural number $n > 0$, and write $\tilde{w} = \Pi_n w$. We denote by $\Phi_t$ the random flow solving (27) and set $w_t = \Phi_t(w)$, $\tilde{w}_t = \Phi_t(\tilde{w})$, $\rho_t = w_t - \tilde{w}_t$. We also use the notation

$$J_t = (D\Phi_t)(w), \qquad \tilde{J}_t = (D\Phi_t)(\tilde{w}), \qquad J_{\rho,t} = J_t - \tilde{J}_t.$$

Since the derivatives of $\phi$ are bounded, the expression inside the limit in (30) is bounded by

$$C \sup_{w \in \mathcal{H}} e^{-\eta \|w\|^2} (\mathbf{E}\|\rho_t\| + \sqrt{\mathbf{E}\|\rho_t\|^2 \mathbf{E}\|J_t\|^2} + \mathbf{E}\|J_{\rho,t}\|).$$

The claim then follows immediately from Theorem A.3 and from the a priori bounds of Lemma A.1.

In order to show that the second claim holds, fix a function $\phi \in \mathcal{C}_0^\infty(\mathcal{H})$ which is of the form $\phi = \tilde{\phi} \circ \Pi_n$ for a $\mathcal{C}_0^\infty$ function $\tilde{\phi}$ and some $n > 0$. It is straightforward to check that there exists a constant $C$ (depending on $\tilde{\phi}$) such that

$$\|\mathcal{P}_t \phi - \phi\|_\eta \leq C \sup_{w \in \mathcal{H}} e^{-\eta \|w\|^2}(\mathbf{E}\|\Pi_n w_t - \Pi_n w\| + \mathbf{E}\|\Pi_n J_t - \Pi_n\|)$$

$$\stackrel{\text{def}}{=} C \sup_{w \in \mathcal{H}} e^{-\eta \|w\|^2} (G_1(t) + G_2(t)).$$

Since $n$ is fixed, both terms are relatively easy to control in the limit $t \to 0$.

Let us first bound $G_1(t)$. It follows from the variation of constants formula (or the mild formulation of a solution) and (37) from the Appendix that

$$G_1(t) \leq \|(1 - \Pi_n e^{\nu \Delta t})\|w\| + \mathbf{E}\left\|\int_0^t \Pi_n e^{\nu \Delta (t-s)} B(\mathcal{K} w_s, w_s)\,ds\right\|$$

$$\leq (1 - e^{-\nu n^2 t})\|w\| + Cn^3 \int_0^t \mathbf{E}\|\Pi_n B(\mathcal{K} w_s, w_s)\|_{-3}\,ds$$

$$\leq (1 - e^{-\nu n^2 t})\|w\| + Cn^3 \int_0^t \mathbf{E}\|w_s\|^2\,ds.$$

Since $n$ is fixed, it is obvious that the first term converges to 0 as $t \to 0$. By (41), $\mathbf{E}\|w_s\|^2$ is uniformly bounded in time by $C \exp(\eta\|w\|^2)$. Hence the second term is bounded by $C \exp(\eta\|w\|^2) t$ and thus converges to 0 as $t \to 0$.

The term $G_2(s)$ is bounded in much the same way. Again it follows from the variation of constants formula that

$$\Pi_n J_t \xi = \Pi_n e^{\nu \Delta t} \xi + \int_0^t e^{\nu \Delta (t-s)} \Pi_n (B(\mathcal{K} J_s \xi, w_s) + B(\mathcal{K} w_s, J_s \xi))\,ds.$$

It follows from (37) that one has the almost sure bound

$$\|\Pi_n J_t - \Pi_n\| \leq 1 - e^{-\nu n^2 t} + Cn^3 \int_0^t \|w_s\|\|J_s\|\,ds.$$

Taking expectations, the needed bound showing that $G_2(t) \to 0$ as $t \to 0$ follows from Lemma A.1 and the same reasoning as used for $G_1(t)$. $\square$

Since the semigroup $\mathcal{P}_t$ is strongly continuous on $\mathcal{B}_\eta$, it has an infinitesimal generator $\mathcal{L}$. Itô's formula allows us to show that $\mathcal{L}$ is an extension of some concrete second-order differential operator:

LEMMA 5.11. *Let $\mathcal{L}$ be the generator of $\mathcal{P}_t$ on $\mathcal{B}_\eta$ and let $\phi \in \mathcal{B}_\eta$ be of the form $\phi(w) = \tilde{\phi} \circ \Pi_n$ for some $n$ and some function $\tilde{\phi} \in \mathcal{C}_0^\infty(\mathbf{R}^n)$. Then $\phi \in \mathcal{D}(\mathcal{L})$ and*

$$
\begin{aligned}
(\mathcal{L}\phi)(w) = {}& \nu \langle \Delta D\phi(w), w \rangle - \langle B(\mathcal{K}w, D\phi(w)), w \rangle \\
& + \langle \bar{f}, D\phi(w) \rangle + \tfrac{1}{2} \operatorname{tr}(QD^2\phi(w)),
\end{aligned}
\tag{31}
$$

*for every $w \in \mathcal{H}$.*

PROOF. Fix a function $\phi$ as in the statement of the lemma. Note first that $D\phi(w) \in \mathcal{D}(\Delta)$ so that (31) does indeed make sense for every $w \in \mathcal{H}$.

One has

$$\Pi_n w_t = \nu \int_0^t \Delta \Pi_n w_s \, ds + \int_0^t \Pi_n B(\mathcal{K}w_s, w_s) \, ds + QW(t),$$

so that Itô's formula immediately implies that

$$\mathcal{P}_t \phi(w) - \phi(w) = \int_0^t \mathcal{P}_s \mathcal{L}\phi(w) \, ds, \tag{32}$$

where $\mathcal{L}\phi$ is given by (31). Let us show that $\mathcal{L}\phi \in \mathcal{B}_\eta$. The only term in (31) for which this is not immediate is the one involving the nonlinearity $B$. Since $D\phi(w) = D\phi(\Pi_m w)$ for $m \geq n$, one has the bound

$$
\begin{aligned}
|\langle B(\mathcal{K}w, D\phi(w)), w\rangle &- \langle B(\mathcal{K}\Pi_m w, D\phi(\Pi_m w)), \Pi_m w\rangle| \\
&\leq |\langle B(\mathcal{K}w - \mathcal{K}\Pi_m w, D\phi(w)), w\rangle| + |\langle B(\mathcal{K}\Pi_m w, D\phi(w)), w - \Pi_m w\rangle| \\
&\leq C\|\mathcal{K}(w - \Pi_m w)\|\|w\|\|D\phi(w)\|_1 + C\|w\|\|D\phi(w)\|_3\|w - \Pi_m w\|_{-1} \\
&\leq \frac{C}{n}\|w\|^2,
\end{aligned}
$$

and similarly for its derivative. The penultimate inequality in this equation is obtained by making use of the bound $\|B(\mathcal{K}w, \tilde{w})\|_1 \leq C\|w\|\|\tilde{w}\|_3$. The result then follows from (32) and the fact that $\mathcal{P}_t$ is strongly continuous. $\square$

5.4. *Convergence of structure functions.* In this section, we show that if $\phi\colon H^1 \to \mathbf{R}$ is a smooth function with at most polynomial growth, then there exist constants $C$, $\eta$ and $\gamma$ (with only $C$ depending on $\phi$) such that

$$\text{(33)} \qquad \left|(\mathcal{P}_t\phi)(w) - \int_{H^1} \phi(w)\,\mu_\star(dw)\right| \leq C e^{\eta\|w\|^2 - \gamma t}.$$

In particular, since $w \in H^1$ implies that $v \in H^2 \subset \mathcal{C}(\mathbf{T}^2, \mathbf{R}^2)$, polynomials of point evaluations of the velocity field fall into this class of observables.

It follows from the results of the previous section that (33) is an immediate consequence of the following result:

PROPOSITION 5.12. *Let $N > 0$ and let $\phi\colon H^1 \to \mathbf{R}$ be a smooth function with*

$$\|\phi\|_N = \sup_{w \in H^1} \frac{|\phi(w)| + |D\phi(w)|}{1 + \|w\|_1^N} < \infty.$$

*Then, for every $t > 0$ and every $\eta > 0$ one has $\mathcal{P}_t\phi \in \tilde{\mathcal{B}}_\eta$. In particular there exist constants $C_{N,t}$ such that $\|\mathcal{P}_t\phi\|_\eta \leq C_{N,t}\|\phi\|_N$.*

PROOF. Fix arbitrary values for $t > 0$ and $\eta > 0$. Let $w \in \mathcal{H}$ and let $w_t$ denote the solution to (SNS) starting at $w$. One then has

$$|\mathcal{P}_t\phi(w)| \leq \|\phi\|_N \mathbf{E}(1 + \|w_t\|_1^N) \leq C\exp(\eta\|w\|^2)\|\phi\|_N,$$

where the second inequality follows from (41). One furthermore has, for an arbitrary vector $\xi \in \mathcal{H}$,

$$|D\mathcal{P}_t\phi(w)\xi| = |\mathbf{E}D\phi(w_t)J_{0,t}\xi| \leq \|\phi\|_N (\mathbf{E}(1 + \|w_t\|_1)^{2N}\mathbf{E}\|J_{0,t}\xi\|_1^2)^{1/2}$$
$$\leq C\exp(\eta\|w\|^2)\|\phi\|_N,$$

where the last bound was obtained by combining (41), (44) and (40). The claim follows immediately from these two estimates. □

5.5. *Regular dependence on the parameters.* In this section, we present one possible application of the results obtained in this article. It was shown in [21] that, for a large class of parameters $\nu$, $Q$ and $\bar{f}$, (SNS) has a unique invariant measure $\mu_\star$. One question which was not addressed was the nature of the dependence of $\mu_\star$ on these parameters. The results obtained in this article enable us to give a relatively simple argument that shows that $\mu_\star$ depends in a continuous way on all the parameters involved. In [32], Majda and Wang proved that in the setting where the dissipation dominates the dynamics, the invariant measure depends continuously on the viscosity. This is a reflection of the fact in this context the random attractor consists of a single point (see [32, 36, 38]). Hence the continuous dependence of the invariant measure follows from the continuous dependence of the random attractor. This can be found in [45] in an abstract setting and [32] in this setting. In the setting of this article, the random attractor is not necessarily a single point, hence results for the attractor do not translate to results for the invariant measure. Nonetheless we show that the long-term statistics of the equations with nearby parameters are near to each other. In particular, our results hold even when the viscosity is not large relative to the typical scale of the energy of the forcing.

In order to keep the notation at a bearable level, we introduce the parameter space $\Lambda = \mathbf{R}_+ \times \ell_+^2 \times \mathcal{H}$ and we denote its elements by

$$\alpha = (\nu, Q, \bar{f}).$$

We equip $\Lambda$ with the natural distance given by

$$d(\alpha, \tilde{\alpha})^2 = |\nu - \tilde{\nu}|^2 + \|Q - \tilde{Q}\|^2 + \|\bar{f} - \tilde{\bar{f}}\|^2.$$

We denote by $\Lambda_0$ the subset of $\Lambda$ that satisfies Assumption 1. For every $\alpha \in \Lambda_0$, we denote by $\mu_\star^\alpha$ the unique invariant measure for (SNS) with parameters $\alpha$ and by $\mathcal{P}_t^\alpha$ the corresponding semigroup. For $\tilde{\alpha} \in \Lambda$, $\mu_\star^{\tilde{\alpha}}$ will simply denote any probability measure invariant, not necessarily unique, for $(\mathcal{P}_t^{\tilde{\alpha}})^*$. One then has the following regularity result:

THEOREM 5.13. *For every $\alpha \in \Lambda_0$, there exist $\eta > 0$, $\varepsilon > 0$ and $C_\alpha > 0$ such that*

$$d_\eta(\mu_\star^\alpha, \mu_\star^{\tilde{\alpha}}) \leq C_\alpha \, d(\alpha, \tilde{\alpha}),$$

*for every $\tilde{\alpha} \in \Lambda$ with $d(\alpha, \tilde{\alpha}) \leq \varepsilon$.*

REMARK 5.14. Going carefully through the proofs of the results in this article and keeping track of the dependence of all a priori estimates on the parameters, we believe that one can show that it is possible to choose for $\eta$, $\varepsilon$ and $C_\alpha$ continuous functions of $\alpha$. The main obstacle to this program is to recover the bounds of [39] under weaker assumptions on $Q$.

REMARK 5.15. Even though $\Lambda_0$ is dense in $\Lambda$, this result does not allow to conclude anything about the set of invariant measures for $\alpha \notin \Lambda_0$. One would expect that there exist values of $\alpha$ such that (SNS) with parameters $\alpha$ has more than one invariant measure. This would then necessarily imply that $C_\beta \gtrsim 1/d(\alpha, \beta)$ for $\beta \in \Lambda_0$ close to $\alpha$.

Theorem 5.13 is the result of the following meta theorem. Given two Markov semigroups, if one is uniformly ergodic and the other is close to the first on $O(1)$ time intervals, then any invariant measure of the second is close to the unique invariant measure of the first. Theorem 1.4 gives the needed ergodicity for $\alpha \in \Lambda_0$. The closeness of the time $t$ transition densities is given by Corollary 5.17 below. It follows from the following bound on the difference between solutions to (SNS) with different sets of parameters:

PROPOSITION 5.16. *Let $w_0 \in \mathcal{H}$ and, for any two sets of parameters $\alpha$ and $\tilde{\alpha}$, let us denote by $w_t$ the solution to* (SNS) *starting at $w_0$ with parameters $\alpha$ and by $\tilde{w}_t$ the solution starting at $w_0$ with parameters $\tilde{\alpha}$.*

*Then, for every $\alpha \in \Lambda$, there exist $\eta_0 > 0$ and $\varepsilon > 0$ such that, for every $\eta \leq \eta_0$ there exist $\gamma > 0$, and $C > 0$ so that*

$$\mathbf{E}\|w_t - \tilde{w}_t\|^2 \leq Ce^{\gamma t + \eta \|w_0\|^2} d(\alpha, \tilde{\alpha})^2,$$

*for every $\tilde{\alpha} \in \Lambda$ with $d(\alpha, \tilde{\alpha}) \leq \varepsilon$.*

We now use this result to prove the needed estimate on the closeness of the time $t$ dynamics.

COROLLARY 5.17. *For any $\alpha \in \Lambda$ there exists an $\eta_0 > 0$ so that for any $\eta \leq \eta_0$ there exist $\gamma > 0$, $\epsilon > 0$, $t_0 > 0$ and $C > 0$ so that one has*

$$d_\eta((\mathcal{P}_t^\alpha)^*\mu, (\mathcal{P}_t^{\tilde{\alpha}})^*\mu) \leq Ce^{\gamma t} d(\alpha, \tilde{\alpha}) \int_\mathcal{H} e^{\eta\|w\|^2} \mu(dw)$$

*for any measure $\mu$ on $\mathcal{H}$, $t \geq t_0$ and $\tilde{\alpha} \in \Lambda$ with $d(\alpha, \tilde{\alpha}) < \epsilon$.*

For brevity in the sequel, we will simply write $\mathcal{P}_t^{\alpha *}$ for $(\mathcal{P}_t^\alpha)^*$.

PROOF OF COROLLARY 5.17. First note that, for every pair $(w, \tilde{w})$ in $\mathcal{H}$ and for every $\eta > 0$, one has the upper bound

(34) $$d_\eta(w, \tilde{w}) \leq \|w - \tilde{w}\|(e^{\eta\|w\|^2} + e^{\eta\|\tilde{w}\|^2}).$$

Fix now $\alpha > 0$, let $\varepsilon$ be as given by Proposition 5.16, and choose an arbitrary $\tilde{\alpha} \in \Lambda$ with $d(\alpha, \tilde{\alpha}) \leq \varepsilon$. Using the notation of Proposition 5.16, we have for $\eta$ sufficiently small

$$d_\eta(\mathcal{P}_t^{\alpha *}\delta_{w_0}, \mathcal{P}_t^{\tilde{\alpha} *}\delta_{w_0}) \leq \mathbf{E}\,d_\eta(w_t, \tilde{w}_t) \leq (\mathbf{E}\|w_t - \tilde{w}_t\|^2 \mathbf{E}(e^{2\eta\|w_t\|^2} + e^{2\eta\|\tilde{w}_t\|^2}))^{1/2}$$

$$\leq C\,d(\alpha, \tilde{\alpha})\exp\left(\gamma t + \frac{\eta}{2}\|w_0\|^2 + \eta e^{((-(\nu-\varepsilon)t)/2)}\|w_0\|^2\right).$$

This shows that there exist constants $t_0$, $\gamma$ and $C$ such that

$$d_\eta(\mathcal{P}_t^{\alpha *}\mu, \mathcal{P}_t^{\tilde{\alpha} *}\mu) \leq C\,d(\alpha, \tilde{\alpha})e^{\gamma t} \int_\mathcal{H} e^{\eta\|w\|^2} \mu(dw),$$

for every $t \geq t_0$. By Remark A.2 we can choose the constants uniform over all $\tilde{\alpha}$ with $d(\alpha, \tilde{\alpha}) \leq \epsilon$. □

With Corollary 5.17 in hand, we return to the proof of Theorem 5.13.

PROOF OF THEOREM 5.13. We know from Theorem 5.5 that there exists $t_1$ such that

$$d_\eta(\mathcal{P}_t^{\alpha*}\mu, \mathcal{P}_t^{\alpha*}\nu) \leq \tfrac{1}{2} d_\eta(\mu, \nu),$$

for every $t \geq t_1$. Let $t_0$ be as in Corollary 5.17. Choosing $t = \max\{t_0, t_1\}$, we have

$$d_\eta(\mu_\star^\alpha, \mu_\star^{\tilde{\alpha}}) = d_\eta(\mathcal{P}_t^\alpha \mu_\star^\alpha, \mathcal{P}_t^{\tilde{\alpha}} \mu_\star^{\tilde{\alpha}}) \leq d_\eta(\mathcal{P}_t^{\alpha*}\mu_\star^\alpha, \mathcal{P}_t^{\alpha*}\mu_\star^{\tilde{\alpha}}) + d_\eta(\mathcal{P}_t^{\alpha*}\mu_\star^{\tilde{\alpha}}, \mathcal{P}_t^{\tilde{\alpha}*}\mu_\star^{\tilde{\alpha}*})$$
$$\leq \tfrac{1}{2} d(\mu_\star^\alpha, \mu_\star^{\tilde{\alpha}}) + d(\alpha, \tilde{\alpha}) e^{\gamma t} \int_{\mathcal{H}} e^{\eta \|w\|^2} \mu_\star^{\tilde{\alpha}}(dw).$$

Notice that in (28) the constants on the right-hand side of the estimate depend contiguously on the parameters for $\alpha \in \Lambda$. Hence it follows from (28) that, for $\eta$ sufficiently small, $\int_\mathcal{H} e^{\eta \|w\|^2} \mu_\star^{\tilde{\alpha}}(dw)$ is bounded uniformly for all $\tilde{\alpha}$ with $d(\alpha, \tilde{\alpha}) \leq \varepsilon$, so that the claim follows. □

We close this section with the proof of Proposition 5.16, which amounts to the continuous dependence on the parameters in $\Lambda$ of the solution operator of (SNS).

PROOF OF PROPOSITION 5.16. Define $\rho_t = w_t - \tilde{w}_t$, $\delta_\nu = \nu - \tilde{\nu}$, $\delta_f = \bar{f} - \tilde{\bar{f}}$ and $\delta_Q = Q - \tilde{Q}$. One then has

$$d\rho_t = (\nu \Delta \rho_t + \delta_\nu \Delta \tilde{w}_t + B(\mathcal{K} w_t, \rho_t) + B(\mathcal{K} \rho_t, \tilde{w}_t) + \delta_f) \, dt + \delta_Q \, dW.$$

At this point, we introduce the stochastic convolution

$$\Psi_t = \int_0^t e^{\nu \Delta (t-s)} \delta_Q \, dW(s),$$

and we set $\bar{\rho}_t = \rho_t - \Psi_t$. This yields for $\bar{\rho}$

$$\tfrac{1}{2} \partial_t \|\bar{\rho}_t\|^2 = -\nu \|\bar{\rho}_t\|_1^2 - \delta_\nu \langle \nabla \bar{\rho}_t, \nabla \tilde{w}_t \rangle + \langle B(\mathcal{K} \bar{\rho}_t, \tilde{w}_t), \bar{\rho}_t \rangle$$
$$+ \langle B(\mathcal{K} w_t, \Psi_t), \bar{\rho}_t \rangle + \langle B(\mathcal{K} \Psi_t, \tilde{w}_t), \bar{\rho}_t \rangle + \langle \delta_f, \bar{\rho}_t \rangle.$$

Fix now $\eta > 0$. Making use of (36), we see that there exists a universal constant $C$ such that

$$\partial_t \|\bar{\rho}_t\|^2 \leq -\nu \|\bar{\rho}_t\|_1^2 + \frac{\delta_\nu^2}{\nu} \|\tilde{w}_t\|_1^2 + C \|\tilde{w}_t\|_1 \|\bar{\rho}_t\|_{1/2} \|\bar{\rho}_t\|$$
$$+ \frac{\eta \nu}{2} (\|w_t\|_1^2 + \|\tilde{w}_t\|_1^2) \|\bar{\rho}_t\|^2 + \frac{C}{\eta \nu} \|\Psi_t\|_1^2 + \langle \delta_f, \bar{\rho}_t \rangle.$$

Note now that it follows from Hölder and Young's inequalities that there exists a universal constant $C'$ such that

$$C\|\tilde{w}_t\|_1 \|\bar{\rho}_t\|_{1/2} \|\bar{\rho}_t\| \le \nu \|\bar{\rho}_t\|_1^2 + \frac{\eta\nu}{2}\|\tilde{w}_t\|_1^2 \|\bar{\rho}_t\|^2 + \frac{C'}{\eta^2\nu^3}\|\bar{\rho}_t\|^2.$$

Combining these bounds yields

$$\partial_t \|\bar{\rho}_t\|^2 \le \left(1 + \frac{C'}{\eta^2\nu^3} + \eta\nu(\|w_t\|_1^2 + \|\tilde{w}_t\|_1^2)\right)\|\bar{\rho}_t\|^2$$
$$+ \|\delta_f\|^2 + \frac{C}{\eta\nu}\|\Psi_t\|_1^2 + \frac{\delta_\nu^2}{\nu}\|\tilde{w}_t\|_1^2.$$

We can now apply Gronwall's inequality to get the bound

$$\|\bar{\rho}_t\|^2 \le \exp\left(\left(1 + \frac{C'}{\eta^2\nu^3}\right)t + \eta\nu\int_0^t (\|w_s\|_1^2 + \|\tilde{w}_s\|_1^2)\,ds\right)$$
$$\times \left(\|\delta_f\|^2 t + \frac{C}{\eta\nu}\int_0^t \|\Psi_s\|_1^2\,ds + \frac{\delta_\nu^2}{\nu}\int_0^t \|\tilde{w}_s\|_1^2\,ds\right).$$

Using the bound $x \le a^{-1}e^{ax}$, applying Cauchy–Schwarz and using the fact that there exists a universal constant $C$ such that, for every Gaussian random variable taking values in a separable Hilbert space, one has

$$\mathbf{E}\|X\|^4 \le C(\mathbf{E}\|X\|^2)^2,$$

we eventually get that there exist constants $C$ and $\gamma$ depending continuously on $\eta$ and on the parameters $\alpha$ and $\tilde{\alpha}$ such that, for every $\eta$ sufficiently small, one has the bound

$$\mathbf{E}\|\bar{\rho}_t\|^2 \le Ce^{\gamma t + \eta\|w_0\|^2}\left(\delta_\nu^2 + \|\delta_f\|^2 + \int_0^t \mathbf{E}\|\Psi_s\|_1^2\,ds\right).$$

The claim then follows immediately from the fact that

$$\mathbf{E}\|\Psi_t\|_1^2 \le \frac{\|\delta_Q\|^2}{2\nu},$$

for every $t \ge 0$. $\square$

**6. Discussion.** We have proven a spectral gap in a Wasserstein distance for a class of Markov processes satisfying a gradient estimate and a weak (topological) irreducibility assumption. Measuring convergence in a Wasserstein metric allows one to incorporate information about the pathwise contraction properties of the system. When the system is completely pathwise contracting, the story is relatively straightforward; see [36, 38] or [25] for the finite-dimensional setting. However, when the system is not pathwise contracting one must introduce a change of measure to make it contracting.

This was one of the central ideas used in [14, 37]. The term in the gradient estimate which does not have a derivative reflects the probabilistic cost of this change of measure while the term with a gradient but a coefficient less than 1 reflects the contraction property obtained via the change of measure.

When the gradient estimate is not uniform, the existence of a Lyapunov function is required. The convergence is then measured in a Wasserstein distance weighted by the Lyapunov function. In this "Harris-like" setting, the contraction properties of the system arise from two sources. Points close to the center of the phase space, as measured by the value of the Lyapunov function, contract due to the combination of deterministic contraction and probabilistic mixing captured by the gradient estimate. Points far out in the space move closer to each other in the distance weighted by the Lyapunov function simply because the linear instability of the flow is compensated by the decrease of the values of the Lyapunov function as the solution moves points toward the center of the phase space.

While we have applied our general theory to the single example of the stochastic Navier–Stokes equations with degenerate forcing, we believe that these results will be useful in many contexts. The gradient estimate allows to capture the combination of mixing due to the presence of noise and due to the contractive nature of the dynamic in one simple estimate. In the context of degenerately forced dissipative SPDEs, control of the gradient term on the right-hand side of Assumption 5 combines an argument strongly inspired by the probabilistic proofs of Hörmander's theorem [24] based on Malliavin's calculus [33, 41, 48], together with the infinitesimal equivalent of the Foias–Prodi-type estimate, namely the fact that the linearized flow contracts all but finitely many directions.

This work has its intellectual roots in many papers. In finite dimensions, spectral gaps in weighted total variation norms like (25) have been obtained for some time [40], but these estimates are of course not uniform when (SNS) is approximated by a sequence of finite-dimensional systems (say by spectral Galerkin approximations). In [46], spaces of observables weighted by Lyapunov functions are used to prove the existence of solutions to infinite-dimensional Kolmogorov equations. The convergence of observables dominated by Lyapunov functions was also given in [27, 38] in the "essentially elliptic" case. The results obtained there were, however, far from what is needed to prove a spectral gap. The convergence results are direct descendants of those developed by many authors in, among others, [6, 14, 20, 28, 34, 37, 38, 43]. All of these works make use of a version of the Foias–Prodi-type estimate [18], introduced in the stochastic context in [35]. The later papers also use a coupling construction to prove convergence. In particular, [20, 37, 38] developed a coupling construction to prove exponential convergence. Though in a less explicit way than its predecessors, the present work makes use of both ideas.

## APPENDIX: PRIORI BOUNDS ON THE DYNAMICS

This appendix is devoted to the proof of the technical estimates used throughout the last two sections of this article. The techniques used to derive these estimates are standard. Even though most of these bounds are probably known to the experts in this field, we have not always been able to find references that state them in the form required here. In particular, we need precise bounds on the difference between the solutions (and their Jacobians) for two nearby initial conditions.

We define for $\alpha \in \mathbf{R}$ and for $w$ a smooth function on $[0, 2\pi]^2$ with mean 0 the norm $\|w\|_\alpha$ by

$$\|w\|_\alpha^2 = \sum_{k \in \mathbf{Z}^2 \setminus \{0,0\}} |k|^{2\alpha} w_k^2,$$

where of course $w_k$ denotes the Fourier mode with wavenumber $k$. Define furthermore $(\mathcal{K}w)_k = -iw_k k^\perp / \|k\|^2$. If $v$, $u_1$ and $u_2$ are as $w$ and $u = (u_1, u_2)$, then $B(u,v) = (u \cdot \nabla)v$. Setting $\mathcal{S} = \{s = (s_1, s_2, s_3) \in \mathbf{R}_+^3 : \sum s_i \geq 1, s \neq (1,0,0), (0,1,0), (0,0,1)\}$ and keeping $u$, $v$ and $w$ as above, then the following relations are useful (cf. [7]):

$$\langle B(u,v), w \rangle = -\langle B(u,w), v \rangle \quad \text{if } \nabla \cdot u = 0, \tag{35}$$

$$|\langle B(u,v), w \rangle| \leq C \|u\|_{s_1} \|v\|_{1+s_2} \|w\|_{s_3}, \quad (s_1, s_2, s_3) \in \mathcal{S}, \tag{36}$$

$$\|B(u,v)\|_\alpha \leq C_\alpha \|u\| \|v\| \quad \text{if } \alpha < -2 \text{ and } \nabla \cdot u = 0, \tag{37}$$

$$\|\mathcal{K}v\|_\alpha = \|v\|_{\alpha-1}, \tag{38}$$

$$\|v\|_\beta^2 \leq \varepsilon \|v\|_\alpha^2 + \varepsilon^{-2((\gamma-\beta)/(\beta-\alpha))} \|v\|_\gamma^2 \tag{39}$$

$$\text{if } 0 \leq \alpha < \beta < \gamma \text{ and } \varepsilon > 0.$$

We start with the following set of a priori bounds, most of which were taken from [21] and [39].

LEMMA A.1. *The solution $w_t$ of the 2D stochastic Navier–Stokes equations in the vorticity formulation satisfies the following bounds:*

1. *There exist constants $C, \eta_\star, \gamma > 0$, such that*

$$\mathbf{E} \exp\left( \nu \int_s^t \eta \|w_r\|_1^2 \, dr - \gamma(t-s) \right) \leq C \exp(\eta \|w_0\|^2), \tag{40}$$

   *for every $t \geq s \geq 0$ and for every $\eta \leq \eta_\star$.*
2. *For every $N > 0$, every $t > 0$ and every $\eta > 0$, there exists a constant $C$ such that one has*

$$\mathbf{E}\|w_t\|_1^N \leq C \exp(\eta \|w_0\|^2), \tag{41}$$

   *for every initial condition $w_0 \in \mathcal{H}$.*

3. *There exist constants $\eta_\star > 0$ and $C > 0$ such that for every $t > 0$ and every $\eta \leq \eta_\star$, the bound*

$$\mathbf{E}\exp(\eta\|w_t\|^2) \leq C\exp(\eta e^{((-\nu t)/2)}\|w_0\|^2) \tag{42}$$

*holds.*

4. *For every $\eta > 0$, there exists a constant $C > 0$ such that the Jacobian $J_{0,t}$ satisfies almost surely*

$$\|J_{0,t}\| \leq \exp\left(\eta \int_0^t \|w_s\|_1^2 \, ds + Ct\right), \tag{43}$$

*for every $t > 0$.*

5. *For every $\eta > 0$ and every $T > 0$, there exists a constant $C$ such that*

$$\int_0^t \|J_{0,s}\xi\|_1^2 \, ds \leq C\|\xi\|^2 \exp\left(\eta \int_0^t \|w_s\|_1^2 \, ds\right), \tag{44}$$

*for every $\xi \in \mathcal{H}$ and every $t \in [0, T]$.*

6. *For every $\eta > 0$ there exists a constant $C$ such that*

$$\|J_{0,t}\xi\|_1^2 \leq C\|\xi\|^2 \exp\left(\eta \int_0^t \|w_s\|_1^2 \, ds + Ct\right), \tag{45}$$

*almost surely, for every $t > 0$ and for every $\xi \in \mathcal{H}$.*

7. *For every $\eta > 0$ and every $p > 0$, there exists $C > 0$ such that the Hessian $K_{0,t}$ satisfies*

$$\mathbf{E}\|K_{0,t}\|^p \leq C\exp(\eta\|w_0\|^2), \tag{46}$$

*for every $t \in [0, 1]$.*

REMARK A.2. It is straightforward to verify that if one fixes a $K_1 > 0$ and $K_2 > 0$, the constants $C$, $\eta_\star$ and $\gamma$ from the statements in Lemma A.1 can be chosen uniformly over all $\nu > K_1$ and $\|Q\|, \|\bar{f}\| \leq K_2$.

PROOF OF LEMMA A.1. Points 1, 4 and 7 are taken from Lemma 4.10 in [21]. Point 2 follows from Lemma A.4 in [39] and point 6 follows from Lemma B.1 in [39]. Point 3 follows immediately from (28).

It remains to show Point 5. It follows from the linearization of the Navier–Stokes equations that

$$\|J_{0,t}\xi\|^2 - \|\xi\|^2 = -2\nu \int_0^t \|J_{0,s}\xi\|_1^2 \, ds + \int_0^t \langle J_{0,s}\xi, B(\mathcal{K}J_{0,s}\xi, w_s)\rangle \, ds.$$

Using (36), this in turn implies that

$$\int_0^t \|J_{0,s}\xi\|_1^2 \, ds \leq \frac{\|\xi\|^2}{2\nu} + \frac{1}{2\nu}\int_0^t \|w_s\|_1 \|J_{0,s}\xi\| \|J_{0,s}\xi\|_1 \, ds$$

$$\leq \frac{\|\xi\|^2}{2\nu} + \frac{1}{8\nu^2}\int_0^t \|w_s\|_1^2 \|J_{0,s}\xi\|^2 \, ds + \frac{1}{2}\int_0^t \|J_{0,s}\xi\|_1^2 \, ds.$$

It thus follows from (43) that

$$\int_0^t \|J_{0,s}\xi\|_1^2 \, ds \leq \frac{\|\xi\|^2}{\nu} + C\|\xi\|^2 \exp\left(\eta \int_0^t \|w_s\|_1^2 \, ds + Ct\right) \int_0^t \|w_s\|_1^2 \, ds,$$

and the result follows immediately. □

In the remainder of this section, we use the following notation, which is the same as in the proof of Theorem 5.10. We fix an element $w \in \mathcal{H}$ and a natural number $n > 0$. We denote by $\Pi_n$ the orthogonal projection in $\mathcal{H}$ onto the Fourier modes with $|k| \leq n$ and we write $\tilde{w} = \Pi_n w$. We denote by $\Phi_t$ the random flow solving (27) and set $w_t = \Phi_t(w)$, $\tilde{w}_t = \Phi_t(\tilde{w})$, $\rho_t = w_t - \tilde{w}_t$. We also use the notation

$$J_t = (D\Phi_t)(w), \qquad \tilde{J}_t = (D\Phi_t)(\tilde{w}), \qquad J_{\rho,t} = J_t - \tilde{J}_t.$$

The aim of this section is to show that, given $t > 0$ and provided $n$ is large enough, it is possible to make $\rho_t$ and $J_{\rho,t}$ arbitrarily small. More precisely, the main result of this section is:

THEOREM A.3. *For every $\gamma > 0$, every $T > 0$ and every $\eta > 0$ there exists $n > 0$ such that*

$$\mathbf{E}\|\rho_T\|^2 \leq \gamma \exp(\eta\|w\|^2), \qquad \mathbf{E}\|J_{\rho,T}\|^2 \leq \gamma \exp(\eta\|w\|^2),$$

*for every $w \in \mathcal{H}$.*

We define the family of increasing stochastic processes $F_\eta^p(t)$ by

$$F_\eta^p(t) = \exp\left(2\eta \int_0^t (\|w_s\|_1^2 + \|\tilde{w}_s\|_1^2) \, ds\right)\left(1 + \sup_{s \in [0,t]} (\|w_s\| + \|\tilde{w}_s\|)^p\right).$$

Note that one has the following result, the proof of which is a trivial application of the a priori bounds from Lemma A.1:

LEMMA A.4. *For every $\eta > 0$, every $t > 0$ and every $p > 0$ there exist $\eta_0 > 0$ and $C$ such that*

$$\mathbf{E}(F_\zeta^p(t)) \leq C \exp(\eta\|w\|^2),$$

*uniformly for every $n \geq 0$, every $w \in \mathcal{H}$ and every $\zeta \in [0, \eta_0]$.*

PROOF OF THEOREM A.3. We fix a terminal time $T > 0$ and start with the bound for $\|\rho_T\|$, which is almost identical to the proof of [21], Lemma 4.17. Note first that $\rho$ solves the equation

$$\partial_t \rho_t = \nu \Delta \rho_t + \tilde{B}(\rho_t, w_t + \tilde{w}_t),$$

where we set $\tilde{B}(w, \tilde{w}) = B(\mathcal{K}w, \tilde{w}) + B(\mathcal{K}\tilde{w}, w)$. Define $\rho_t^\ell = \Pi_n \rho_t$ and $\rho_t^h = \rho_t - \rho_t^\ell$, so that

$$\partial_t \|\rho_t^\ell\|^2 = -2\nu \|\rho_t^\ell\|_1^2 + \langle B(\mathcal{K}\rho_t^\ell, w_t + \tilde{w}_t), \rho_t^\ell \rangle$$
$$- \langle B(\mathcal{K}\rho_t^h, \rho_t^\ell), w_t + \tilde{w}_t \rangle - \langle B(\mathcal{K}w_t + \mathcal{K}\tilde{w}_t, \rho_t^\ell), \rho_t \rangle,$$
$$\partial_t \|\rho_t^h\|^2 = -\nu \|\rho_t^h\|_1^2 - \langle B(\mathcal{K}\rho_t, \rho_t^h), w_t + \tilde{w}_t \rangle - \langle B(\mathcal{K}w_t + \mathcal{K}\tilde{w}_t, \rho_t^h), \rho_t \rangle.$$

The initial conditions for these equations are given by

$$\rho_0^\ell = 0, \qquad \rho_0^h = \Pi_n w.$$

The equations satisfied by $\rho_t^\ell$ and $\rho_t^h$ are the same as the ones appearing in the proof of [21], Lemma 4.17, so that we get the bounds

$$\|\rho_t^h\|^2 \leq \|w\|^2 \left( e^{-\nu n^2 t} + \frac{C_\eta}{n} F_\eta^1(t) \right),$$
$$\|\rho_t^\ell\|^2 \leq C_\eta \int_0^t \exp\left( \eta \int_s^t \|w_r + \tilde{w}_r\|_1^2 \right) \|w_s + \tilde{w}_s\|_{1/2}^2 \|\rho_s^h\|^2 \, ds$$
$$\leq C_\eta F_\eta^4(t) \int_0^t \|w_s + \tilde{w}_s\|_1 \|\rho_s^h\| \, ds.$$

These bounds are valid for every $\eta > 0$. It follows from the first bound that

$$\int_0^T \|\rho_s^h\|^2 \, ds \leq \frac{C}{n} F_\eta^3(T),$$

so that the second bound yields

(47) $$\sup_{t \in [0,T]} \|\rho_t\|^2 \leq \frac{C_\eta}{\sqrt{n}} F_{2\eta}^6(T).$$

The bound on $\mathbf{E}\|\rho_T\|^2$ then follows from Lemma A.4.

In order to bound $J_{\rho,T}$, note first that $J_{\rho,0} = 0$ and

$$\partial_t J_{\rho,t} = \nu \Delta J_{\rho,t} + \tilde{B}(J_{\rho,t}, w_t + \tilde{w}_t) + \tilde{B}(J_t + \tilde{J}_t, \rho_t).$$

Fix now a tangent vector $\xi \in \mathcal{H}$. It follows from (36) that

$$\partial_t \|J_{\rho,t}\xi\|^2 \leq -2\nu \|J_{\rho,t}\xi\|_1^2 + C\|J_{\rho,t}\xi\|_{1/4} \|J_{\rho,t}\xi\| \|w_t + \tilde{w}_t\|_1$$
$$+ C\|J_{\rho,t}\xi\|_1 \|\rho_t\| \|J_t\xi + \tilde{J}_t\xi\|_{1/4}$$
$$\leq (C_\eta + \eta \|w_t + \tilde{w}_t\|^2) \|J_{\rho,t}\xi\|^2 + \|\rho_t\|^2 \|J_t\xi + \tilde{J}_t\xi\|_{1/4}^2.$$

This bound is valid (with different values for the constant $C_\eta$) for any value of $\eta > 0$. It immediately implies that

$$\|J_{\rho,T}\xi\|^2 \leq F_\eta^0(T) \int_0^T \|\rho_t\|^2 \|J_t\xi + \tilde{J}_t\xi\|_{1/2} \|J_t\xi + \tilde{J}_t\xi\| \, dt$$

$$\le CF_{3\eta}^2(T)\|\xi\| \int_0^T \|\rho_t\| \|J_t\xi + \tilde{J}_t\xi\|_{1/2}\, dt$$

$$\le CF_{4\eta}^2(T)\|\xi\|^{3/2} \int_0^T \|\rho_t\| \|J_t\xi + \tilde{J}_t\xi\|_1^{1/2}\, dt,$$

where we made use of (43). It follows that there exists a constant $C$ such that, for every $\alpha > 0$, one has the bound

$$\|J_{\rho,T}\xi\|^2 \le \left(\frac{1}{\alpha} \int_0^T \|\rho_t\|^2\, dt + \alpha CF_{3\eta}^8(T)\right)\|\xi\|^2 + \alpha \int_0^T (\|J_t\xi\|_1^2 + \|\tilde{J}_t\xi\|_1^2)\, dt.$$

It follows from (44) that

$$\|J_{\rho,T}\|^2 \le \left(\frac{1}{\alpha} \int_0^T \|\rho_t\|^2\, dt + \alpha CF_{3\eta}^8(T)\right),$$

so that the claim follows by combining Lemma A.4 with the bound (47). $\square$

**Acknowledgments.** We would like to thank H. Bray, D. Dolgopyat, G. Friesecke, X.-M. Li, C. Odasso, A.-S. Sznitman, J. York and L.-S. Young for useful discussions and comments. We would like to thank A. Majda and X. Wang for pushing the authors to explicitly include forcing with nonzero mean. We also thank Z. Brzeźniak and an anonymous referee for useful and pertinent comments on the first circulated version of this article.

DEPARTMENT OF MATHEMATICS  
UNIVERSITY OF WARWICK  
COVENTRY CV4 7AL  
UNITED KINGDOM  
E-MAIL: M.Hairer@Warwick.ac.uk

DEPARTMENT OF MATHEMATICS  
DUKE UNIVERSITY  
DURHAM, NORTH CAROLINA 27708  
USA  
E-MAIL: jonm@math.duke.edu